\documentclass[runningheads,a4paper]{llncs}
\usepackage{amsmath}
\usepackage{amssymb}
\usepackage{url}
\usepackage{graphicx}
\usepackage{hyperref}
\usepackage{color}
\usepackage{xcolor}
\usepackage{wrapfig}
\usepackage{siunitx}
\usepackage{array,multirow}
\newcolumntype{B}[1]{>{\centering\arraybackslash}m{#1}}

\newcommand{\keywords}[1]{\par\addvspace\baselineskip
\noindent\keywordname\enspace\ignorespaces#1}

\newcommand{\vc}[1]{\mathbf{#1}}
\newcommand{\mt}[1]{\mathbf{#1}}
\newcommand{\ellp}[1]{\ell_{#1}}

\newcommand{\herm}{*}

\newcommand{\target}{x} 
\newcommand{\vtarget}{\vc{\target}}

\newcommand{\light}{d}
\newcommand{\vlight}{\vc{\light}}

\newcommand{\image}{p} 
\newcommand{\vimage}{\vc{\image}} 

\newcommand{\targetLog}{r}
\newcommand{\vtargetLog}{\vc{\targetLog}}

\newcommand{\lightLog}{l}
\newcommand{\vlightLog}{\vc{\lightLog}}

\newcommand{\imageLog}{s}
\newcommand{\vimageLog}{\vc{\imageLog}}

\newcommand{\coeff}{\alpha} 
\newcommand{\vcoeff}{\boldsymbol{\coeff}} 

\newcommand{\meas}{\vc{y}} 
\newcommand{\Id}{\text{Id}}
\newcommand{\Meas}{\mt{A}}
\newcommand{\Dict}{\mt{\Psi}}
\newcommand{\Fourier}{\mt{F}}

\newcommand{\defgr}{\mathrel{\mathop:\!\!=}}
\begin{document}

\mainmatter 

\title{ Reconstruction Methods in THz Single-pixel Imaging }

\titlerunning{Reconstruction Methods in THz Single-pixel Imaging }

%
\author{Martin Burger, Lea F\"ocke, Lukas Nickel, Peter Jung, Sven Augustin}
\authorrunning{Reconstruction Methods for THz Single-pixel Imaging}

\institute{}

%
%

\toctitle{Lecture Notes in Computer Science}
\tocauthor{Authors' Instructions}
\maketitle

\begin{abstract}
The aim of this paper is to discuss some advanced aspects of image reconstruction in single-pixel cameras, focusing in particular on detectors in the THz regime.
We discuss the reconstruction problem from a computational imaging perspective and provide a comparison of the effects of several state-of-the art regularization techniques.

Moreover, we focus on some advanced aspects arising in practice with THz cameras, which lead to nonlinear reconstruction problems: the calibration of the beam reminiscent of the Retinex problem in imaging and phase recovery problems.
Finally we provide an outlook to future challenges in the area.
\end{abstract}

\section{Introduction}
Imaging science has been a strongly evolving field in the last century, with a lot of interesting developments concerning devices, measurement strategies and computational approaches to obtain high-quality images.
A current focus concerns imaging from undersampled data in order to allow novel developments towards dynamic and hyperspectral imaging, where time restrictions forbid to acquire full samplings.
In order to compensate for the sampling either physical models (e.g. concerning motion in dynamic imaging, cf. \cite{burger2018variational}) or a-priori knowledge about the images to be reconstructed are used.
In applications, where one has sufficient freedom to choose the undersampling patterns, the paradigm of compressed sensing is particularly popular.
It is based on minimizing coherence between measurements (cf. \cite{donoho2006compressed}), often achieved by random sampling (cf. \cite{baraniuk2007compressive}).

Single-pixel imaging, or more precisely single-detector imaging, is one of the most interesting developments in compressed sensing (cf. \cite{duarte2008single} \cite{chan2008single} \cite{edgar2018principles} \cite{stokoe2018theory}).
It is based on using a single detector with multiple random masks in order to achieve the desired resolution.
The ideal model is to have the mask realize a grid on the detector region with subpixels of the desired image resolution, which are either open or closed with a certain probability.
The detector integrates the light passing through the open subpixels in the mask.
Note that the image at desired resolution might be acquired by scanning in a deterministic way, in the easiest setting with a single subpixel open at each shot.
However, the light intensity obtained from a single subpixel might not be sufficient to obtain a reasonable signal-to-noise ratio and obviously the mechanical scanning times may strongly exceed the potential times of undersampling with masks having multiple open subpixels.

This idea is particularly relevant in applications where detectors are expensive or difficult to miniaturize such as imaging in the THz range (cf. \cite{chan2008single} \cite{watts2014terahertz} \cite{augustin2015compressed} \cite{augustin2017optically}), which is our main source of motivation.
The random masks are achieved by a spatial light modulator, which may however deviate from the ideal setting in practical applications due to the following effects:
\begin{itemize}

\item {\em Beam calibration:} as in many practical imaging approaches, the lighting beam is not homogeneous and needs to be corrected, a problem reminiscent of the classical Retinex problem (cf. \cite{Land:retinex}).

\item {\em Diffraction:} deviations between the object and image plane may cause the need to consider diffraction and take care of out-of phase effects, which complicates the inversion problem and

\item {\em Motion: } The object to be imaged may move between consecutive shots, i.e. while changing the masks, which leads to well-known motion blur effects in reconstructions.

\end{itemize}
In the following we will discuss such issues in image reconstruction arising in THz single pixel cameras.
We start with the basic modeling of the image formation excluding calibration and further effects in Section 2, where we also discuss several regularization models.
In Section 3 we discuss some computational methods to efficiently solve the reconstruction problem and in particular compare the results of some state-of-the art regularization techniques.
Then we discuss the challenges that arise when applying the single-pixel imaging approach in a practical setup with THz detectors and proceed to calibration models for beams in Section 5, where we provide a relation to the classical Retinex problem and discuss the particular complications arising due to the combination with the single pixel camera.
In Section 6 we discuss the phase reconstruction problem and compare several state-of-the art approaches for such.
Finally we conclude and discuss further challenges in Section 7.

\section{Compressed Sensing and Reconstruction in Single Pixel Cameras}\label{sec:spc}

\begin{figure}[t]
	\centering
	\includegraphics[width=.9\textwidth]{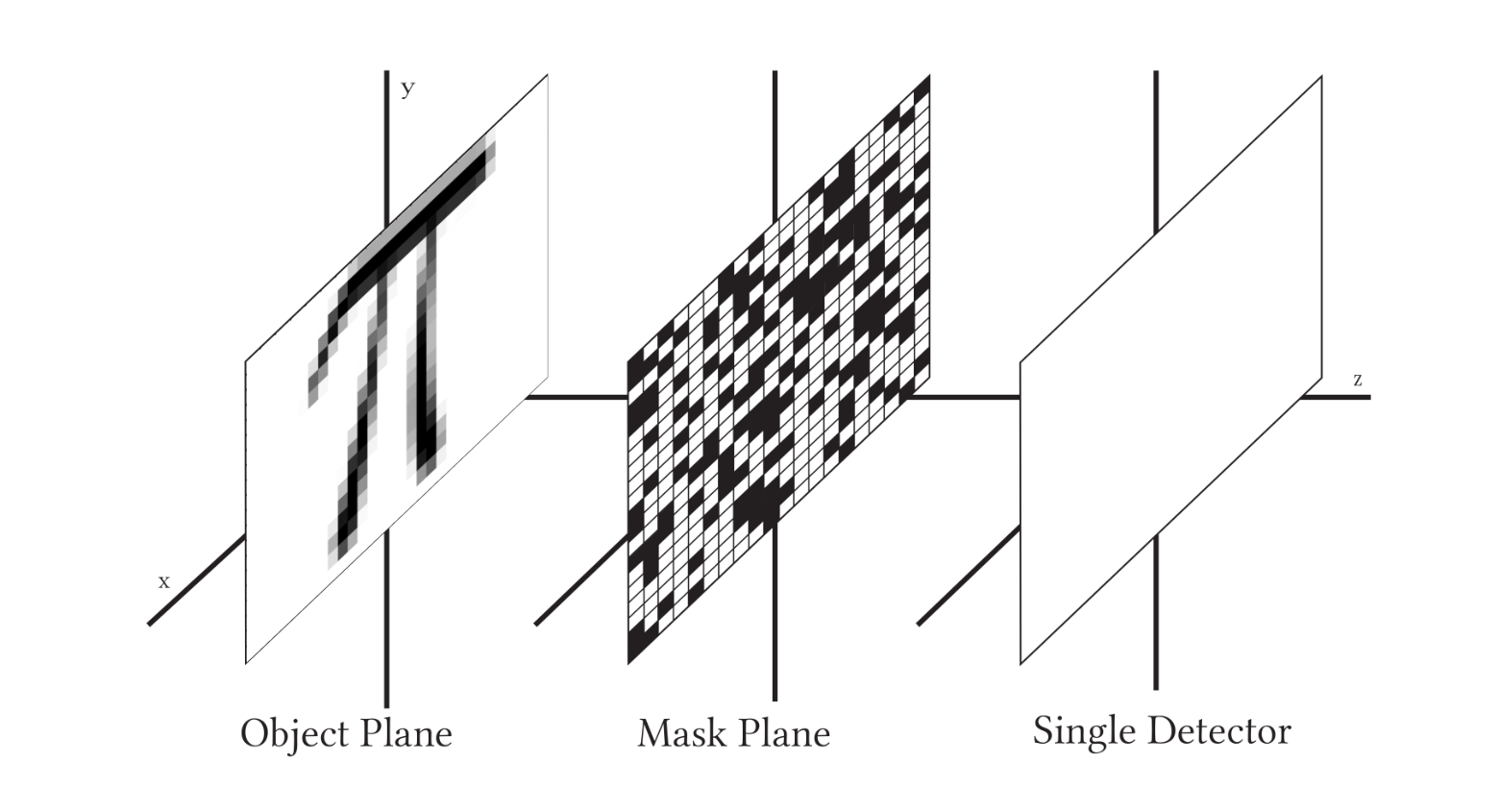}
	\caption{Sketch of the setup of imaging planes in THz single pixel imaging (from \cite{Nickel}).\label{diffractionfigure}}
\end{figure}

In the following we discuss some basic aspects of the image reconstruction in single pixel cameras.
We are interested in reconstructing a two-dimensional image on the subpixel grid of size $d_1\times d_2$, i.e. 
$\vimage\in\mathbb{R}^n$ where $n=d_1\cdot d_2$.
The simplest model of the image formation process is to have each measurement $y_i$, $i=1,\ldots,m$,
as the sum of subpixel values $p_j$ for those $j$ corresponding to the open subpixels in the $i$-th mask.
This means we can write
\begin{equation} \meas = \Meas \vimage, \end{equation}
with a matrix $\Meas \in \mathbb{R}^{m \times n}$ whose entries are only zeros and ones.
Again the non-zero entries in each row of $\Meas$ correspond to the open subpixels of the respective mask (compare Figure \ref{diffractionfigure} for a schematic overview of such a setup).

Choosing $m\geq n$ deterministic masks appropriately one could guarantee that $\Meas$ is invertible and simply solve the linear reconstruction problem.
However, in practice one would like to reconstruct a reasonable image $\vimage$ from a number of measurements $m$ significantly smaller than $n$.
Single pixel cameras are hence developed in the paradigm of compressed sensing and the natural way to realize appropriate measurements is to choose the masks randomly.
In most systems such as the ones we consider here, each entry of $\Meas$ is chosen independently as a binary random variable with fixed expectation.
Combining such approaches with appropriate prior information in terms of sparsity leads to compressed sensing schemes that can be analyzed rigorously (cf. \cite{candes2011probabilistic} \cite{lu2015compressed}).

\subsection{Compressed Sensing Techniques}

A key motivation for compressed sensing comes from the fact that
in many cases $d_1\times d_2$ images
are compressible and can be (approximately) sparsely represented as
$\vimage=\Dict\vcoeff$ with $\vcoeff\in\mathbb{R}^n$ and a particular basis $\Dict\in\mathbb{R}^{n\times n}$, e.g. wavelets, or overcomplete systems $\Dict\in\mathbb{R}^{n\times N}$ with $N > n$ such as shearlets (cf. \cite{guo2007optimally}).
By this we mean that
$$\|\vcoeff\|_{\ellp{0}}=|\{k\,:\,\coeff_k\neq 0\}|, $$
respectively
$$|\{k\,:\,\vert\coeff_k\vert\geq \epsilon\}|, $$
for small $\epsilon > 0$,
is considerably smaller then the ambient dimension $n$.
Ideally, one would thus solve the problem of minimizing $\|\vcoeff\|_{\ellp{0}}$ subject to the linear constraint $\Meas \Dict \vcoeff= \meas$, which is however an NP-hard combinatorial problem.

One of the fundamental results that initiated the field of compressed
sensing (cf. \cite{candes_robust_2006,donoho_compressed_2006,eldar_compressed_2012,foucart_mathematical_2013}) is that under additional assumptions on the map $\Meas\Dict$
the convex relaxation
\begin{equation}
 \min_{\vcoeff} \|\vcoeff\|_{\ellp{1}}\quad\text{s.t.}\quad
 \meas=\Meas\Dict\vcoeff\label{eq:basicCSProblem}
\end{equation}
recovers exactly (in the noiseless setting) the unknown $\vcoeff$ yielding
the correct image $\vimage$.
A common additional assumption is that
the image itself has non-negative pixel intensities, $\vimage\geq 0$
such that problem \eqref{eq:basicCSProblem} is extended to:
\begin{equation}
 \min_{\vcoeff} \|\vcoeff\|_{\ellp{1}}\quad\text{s.t.}\quad
 \meas=\Meas\Dict\vcoeff\quad\text{and}\quad\Dict\vcoeff\geq 0.\label{eq:basicCSProblemWithPos}
\end{equation}
If $\Dict=\Id$ and if the row span of $\Meas$ intersects the positive
orthant, meaning that there exists a vector $\vc{t}$ such that
$\Meas^\herm\vc{t}>0$, the measurement matrix $\Meas$ itself already
assures that the $\ell_1$-norm $\|\vimage\|_{\ellp{1}}$ of a feasible
$\vimage=\vcoeff$ in \eqref{eq:basicCSProblemWithPos} is equal (or
close in the noisy case) to $\|\vimage^0\|_{\ellp{1}}$ where $\vimage^0$  is the
unknown image to recover.
To see this, let us assume
exemplary that we find a vector $\vc{t}$ such that
$\Meas^\herm\vc{t}=\vc{1}$ is the all-one vector and
$\meas=\Meas\vimage^0$.
Then
\begin{equation*}
\begin{split}
\|\vimage\|_{\ellp{1}}-\|\vimage^0\|_{\ellp{1}}
&=
\langle \vc{1},\vimage-\vimage^0\rangle
=\langle \vc{t},\Meas(\vimage-\vimage^0)\rangle\\
&=\langle \vc{t},\Meas\vimage-\meas\rangle
\leq\|\vc{t}\|_{\ellp{2}}\|\Meas\vimage-\meas\|_{\ellp{2}}\\
\end{split}
\end{equation*}
and hence, it is enough to minimize the residual over $\vimage\geq 0$
and replace \eqref{eq:basicCSProblemWithPos} by a simple non-negative
least squares (NNLS) problem:
\begin{equation}
 \min_{\vimage\geq 0} \|\meas-\Meas\vimage\|_{\ellp{2}}^2.\label{eq:basicNNLS}
\end{equation}
Indeed, that under these assumptions $\ell_1$-minimization reduces to
a feasibility problem has observed already in prior work
\cite{Bruckstein2,Wang11}.
In particular the setting of random binary
masks has been investigated in \cite{Slawsky13,kueng:16:nnls} and a considerably (partially-) derandomized result based on orthogonal arrays is discussed in \cite{Jung:frontierts18}.
Although this explains to some extend why
non-negativity and NNLS are very useful in certain imaging problems
this does not easily extends to generic dictionaries $\Dict$.

Hence, coming back to \eqref{eq:basicCSProblemWithPos}, in the case of noisy data knowledge about the expected solution should be included.
Hence one usually rather solves
\begin{equation}
 \min_{\vcoeff} \frac{1}2\|\meas-\Meas\Dict\vcoeff\|_{\ellp{2}}^2+ \lambda \|\vcoeff\|_{\ellp{1}}\quad\text{s.t.} \quad \Dict\vcoeff\geq 0,
\end{equation}
where $\lambda > 0$ is an appropriate regularization parameter.

\subsection{Total Variation Regularization and Related Methods}

While simple $\ell_1$-regularization is a common approach used for the theory of compressed sensing, for image reconstruction this choice as data fitting term has some drawbacks in practice, e.g. in wavelets the reconstructions may suffer from artifacts due to rectangular structures in their constructions.
In more advanced models like shearlets (cf. \cite{kutyniok2012shearlets}) the visual quality of reconstructions is improved, but the computational overhead in such large dictionaries may become prohibitive.
A much more popular approach in practical image reconstructions are total variation based methods such as minimizing
\begin{equation} \label{eq:TVmodel}
 \min_{\vimage} \|\nabla\vimage\|_{\ellp{1}}\quad\text{s.t.}\quad
 \meas=\Meas\vimage
\end{equation}
or the penalty version
\begin{equation}
 \min_{\vimage} \frac{1}2\|\meas-\Meas\vimage\|_{\ellp{2}}^2+\lambda \|\nabla\vimage\|_{\ellp{1}} ,
\end{equation}
potentially with additional non-negativity constraints.
Total variation methods in compressed sensing have been
investigated recently in further detail (cf. \cite{needell2013stable}
\cite{needell2013near} \cite{poon2015role} \cite{Cai:tv15} \cite{krahmer2017total}).

As in wavelet systems, simple approaches in total variation may suffer from rectangular grid artefacts.
This is the case in particular if the straight $\ellp{1}$-norm of $\nabla \vimage$ is used in the regularization (namely $\|\nabla\vimage\|_{\ellp{1,1}}$), which corresponds to an anisotropic total variation.
It is well-known that such approaches promote rectangular structures aligned with the coordinate axis (cf. \cite{caselles2014total}), but destroy round edges.
An improved version is the isotropic version, which considers $\nabla \vimage \in \mathbb{R}^{n \times 2}$ (the rows corresponding to partial derivatives) and computes $\|\nabla\vimage\|_{\ellp{2,1}}$ as total variation.
The isotropic total variation promotes round structures that are visually more appealing in natural images, which can be made precise in a continuum limit.

A remaining issue of total variation regularization in applications is the so-called stair-casing phenomenon, which means that piecewise constant structures are promoted too strongly and thus gradual changes in an image are rather approximated by piecewise constants in a stair-like fashion.
In order to cure such issues infimal convolutions of total variation and higher order functionals are often considered, the most prominent one being the total generalized variation (TGV) model (cf. \cite{bredies_total_2010})
\begin{equation} \label{eq:TGVmodel}
 \min_{\vimage, {\bf w}} \|\nabla\vimage - {\bf w}\|_{\ellp{1}} + \beta \|\nabla {\bf w} \|_{\ellp{1}} \quad\text{s.t.}\quad \meas=\Meas\vimage,
\end{equation}
where $\beta>0$ is a parameter to be chosen appropriately.
Again in practice suitable isotropic variants of the $\ellp{1}$-norm are used, and in most approaches only the symmetrized gradient of the vector field ${\bf w}$ is used instead of the full gradient.

In the penalized form for noisy data used in practice, the models above can be written in the unified form 
$$ \min_{\vimage} \frac{1}2\|\meas-\Meas \vimage\|_{\ellp{2}}^2+ \lambda R({\bf \Phi} \vimage ), $$
 with $R$ being an appropriately chosen seminorm and a suitable matrix ${\bf \Phi}$, e.g. the gradient or ${\bf \Phi} = {\bf \Psi}^{-1}$ in the case of using a basis, or ${\bf \Phi} = {\bf \Psi}^*$ for analysis sparsity.
 A common issue in problems of this form that aim to reduce variance is an increase of bias, which e.g. leads to a loss of contrast and small structures in total variation regularization.
 In order to improve upon this issue iterative regularization (cf. \cite{burger2013guide} \cite{benning2018modern}) can be carried out, in the simplest case by the Bregman iteration, which computes each iteration $\vimage^{k+1}$ as the minimizer of 
$$ \min_{\vimage} \frac{1}2\|\meas-\Meas \vimage\|_{\ellp{2}}^2+ \tilde \lambda (R({\bf \Phi} \vimage ) - {\bf q}^k \cdot \vimage )$$
with subgradient $ {\bf q}^k \in \partial R({\bf \Phi}\vimage^k)$.
In this way in each iteration a suitable distance to the last iterate is penalized instead of the original functional $R$, which is effectively a distance to zero.
The parameter $\tilde \lambda$ is chosen much larger than the optimal $\lambda$ in the above problem, roughly speaking when carrying out $K$ iterations we have $\tilde \lambda = K \lambda.$
 Observing from the optimality condition that $ {\bf q}^k = \mu \Meas^* {\bf r}^k$ with $\mu =\frac{1}{\tilde \lambda}$, the Bregman iteration can be interpreted equivalently as the augmented Lagrangian method for the constrained problem of minimizing $R({\bf \Phi} \vimage )$ subject to $\meas=\Meas \vimage$, i.e. 
$$ \vimage_{k+1} \in \text{arg}\min_{\vimage} \frac{\mu}2\|\meas-\Meas \vimage\|_{\ellp{2}}^2+ R({\bf \Phi} \vimage ) + {\bf r}^k \cdot ( \meas-\Meas \vimage )$$
and 
$$ {\bf r}^{k+1} = {\bf r}^k + \meas-\Meas \vimage^{k+1}.$$

While single pixel cameras are naturally investigated from a compressed sensing point of view, let us comment on some aspects of the problem when viewed as an inverse problem as other image reconstruction tasks in tomography or deblurring.
First of all, we can see some common issues in computational techniques, which are needed due to the indirect relation between the image and the measurements and the fact that the matrix $\Meas$ is related to the discretization of an integral operator.
On the other hand there is a significant difference between the single pixel setup and other image reconstruction problems in the sense that the adjoint operator is not smoothing, i.e. $\Meas^*$ has no particular structure.
Thus, the typical source condition $\Meas^* {\bf{w}} \in {\bf{\Phi}}^* \partial R({\bf{\Phi}} \vimage ) $ for some ${\bf{w}}$ that holds for solutions of the variational problems (cf. \cite{benning2018modern}) does not imply smoothness of the subgradients as in other inverse problems.

\section{Computational Image Reconstruction}

In the following we discuss the effects of different regularization models on the image reconstruction quality.
In order to efficiently compute solutions of the arising variational problems at reasonable image resolution, appropriate schemes to handle the convex but nondifferentiable terms are needed.
It has become a standard approach in computational imaging of such problems to employ first-order splitting methods based on proximal maps that can be computed exactly, we refer to \cite{burger2016first} for an overview.
The key idea is to isolate matrix-vector multiplications and local nonlinearities, for which the proximal map can be computed exactly or numerically efficient.
A standard example is the $\ell_1$-norm, whose $\ell_2$-proximal map 
$$ \text{prox}_{\ell_1}(f) = \text{arg}\min_x \frac{1}2 \Vert x - f\Vert_{\ell_2}^2 + \Vert x \Vert_{\ell_1} $$
is given by soft shrinkage.
A very popular approach in current computational imaging are first-order primal dual methods (cf. \cite{chambolle2011first} \cite{zhang2011unified}) for computing minimizers of problems of the form 
$$ \min_\vimage F(\vimage) + G({\bf L} \vimage) ,$$
with convex functionals $F$ and $G$ and a linear operator ${\bf L}$.
The primal dual approach reformulates the minimization as a saddle point problem
$$ \min_\vimage \max_{\bf q} F(\vimage) + \langle {\bf L} \vimage, {\bf q} \rangle -G^*({\bf q} ) ,$$
with $G^*$ being the convex conjugate of $G$.
Then one iterates primal minimization with respect to $\vimage$ and dual maximization with respect to ${\bf q}$ with additional damping and possible extrapolation steps.
Here we use a toolbox by Dirks \cite{dirks2016flexible} for rapid prototyping with such methods to implement the models of the previous section.

\begin{figure}[t]
	\centering
	\includegraphics[width=.9\textwidth]{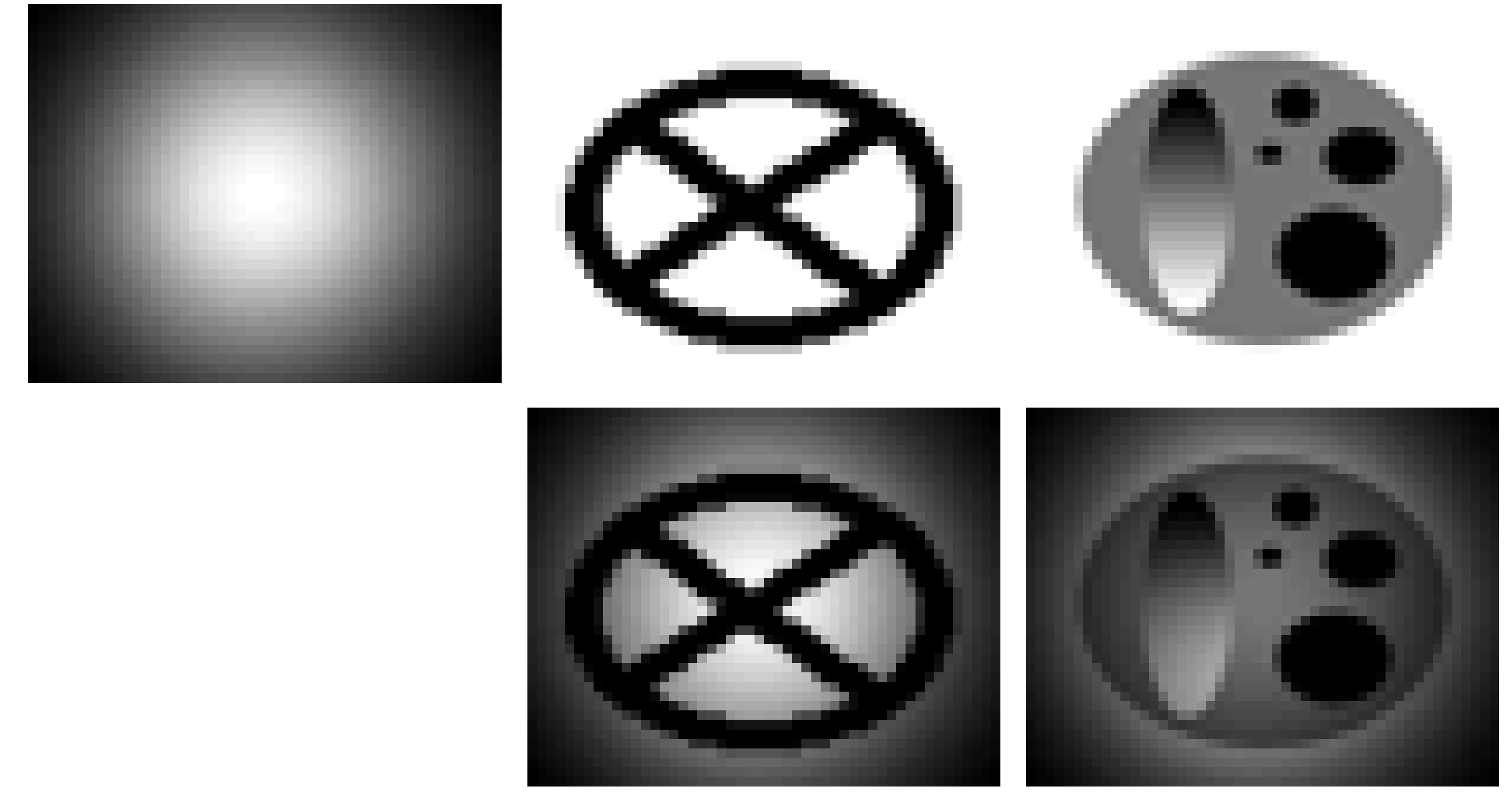}
	\caption{Grayscale images of synthetic data, synthetic illumination beam with values between 0 (black) and 1 (white); top left: Gaussian beam; top center: x phantom ground truth; bottom center: x phantom with applied beam; top right: gradient phantom ground truth; bottom right: gradient phantom with applied beam.}
	\label{im:groundTruth}
\end{figure}

We will use synthetic data to investigate the regularization methods described above.
With this we can fully explore and understand the effects of named methods for different types of images.
Here we consider the multiple regularization operators, namely Tikhonov, $\ell_1$- regularization, total variation (TV) and total generalized variation (TGV) \cite{bredies_total_2010}.
Moreover we compare the approach with a simple non-negative least squares approach that only enforces non-negativity of the image.
Each dataset is composed of a ground truth structural image and a Gaussian kernel beam that is used as an illumination source.
In this work we will explore two different ground truth images, that are illuminated by the shown beam.
Structural images, beam and illuminated structures are shown in Figure \ref{im:groundTruth}.

\begin{figure}[t]
	\centering
	\includegraphics[width=.9\textwidth]{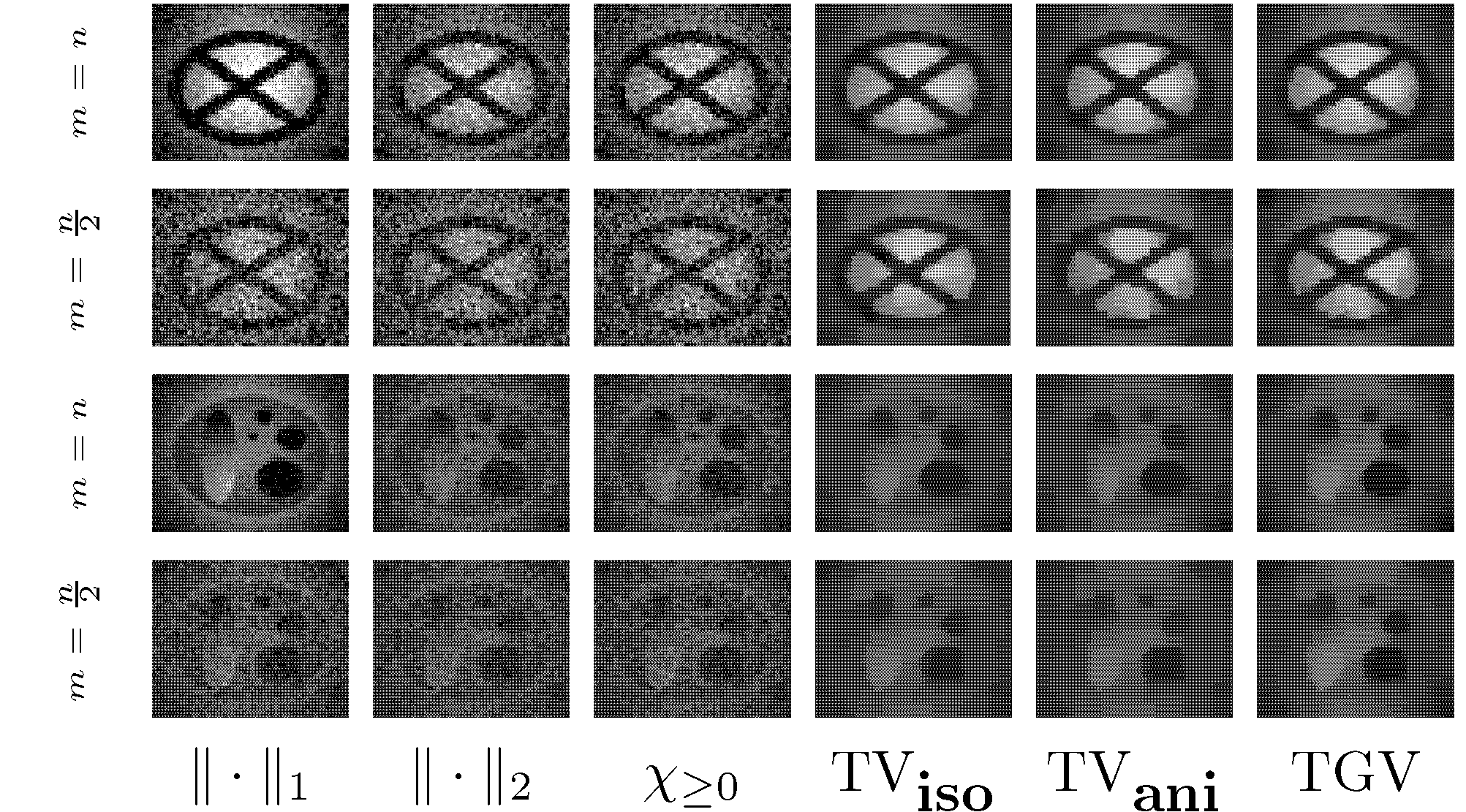}
	\caption{Reconstruction the x and gradient phantom for multiple sampling rates (rows) and reconstruction frameworks (columns).}
	\label{im:compareReconstructions}
\end{figure}

An overview of basic reconstructions in the single detector framework using random binary masks is shown in Figure \ref{im:compareReconstructions}.
Here we use two phantoms and two different sampling ratios.
We display reconstructions with regularization parameters optimized for SSIM and observe the impact of the regularization for both choices of the sampling rate.
In particular we see that for total variation type regularizations there is a strong improvement in visual image quality and that there is hardly any loss when going from full sampling to undersampling.

\section{Challenges in Practical THz Single-Pixel Imaging}

%
%
%
%
 %

\begin{figure}[t]
\centering
 \includegraphics[width=0.4\textwidth]{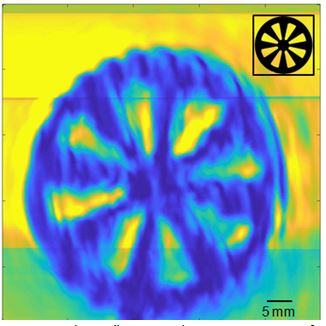}\vspace{-5pt}
 \caption{Mechanically scanned \SI{0.35}{\tera\hertz} image of a metal Siemensstar test target with a diameter of \SI{50}{\milli\meter}.
 	The image acquisition took several hours (\SI{12}{\hour}) depending on the number of steps but still scanning artifacts are very prominent.
The inset shows a photo of the metal Siemensstar target.}
 \label{fig:artifacts}
\end{figure} 

\begin{figure}[t]
\centering
\includegraphics[scale=0.6]{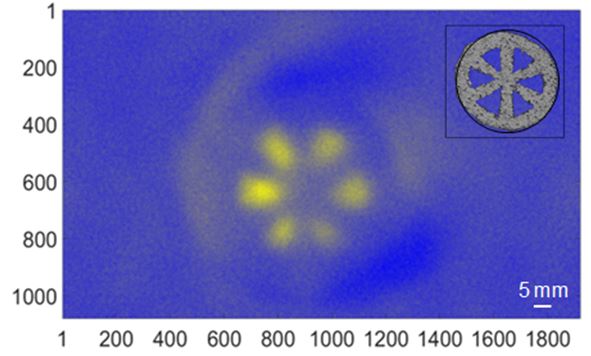}\vspace{-5pt}
\caption{Megapixel \SI{0.35}{\tera\hertz} image using a convolution approach (i.e. a binary circulant matrix).
	Labels show the pixel numbers in both directions, about 10 pixel correspond to $5$mm.
	The spatial resolution in the image is still limited by sub-optimal illumination patterns but the image shows that the resolution is reasonable using the convolution approach, while being very fast (more than one FPS can be achieved).
	The inset shows again a photo of the imaged object.}
\label{fig:convolution}
\end{figure}

If compared to the visible region of the electromagnetic spectrum (VIS) the THz region is quite demanding from a physical point of view.
Due to the large wavelength that is 2-3 orders of magnitude larger than for VIS radiation, the THz region is plagued by coherence effects and diffraction.
This has to be mitigated already on the physical level by using specifically designed optical
elements and techniques.
As a rule of thumb, techniques and methods from the VIS region can be used but they need to be adapted.
Also, lenses and mirrors are used but they are made of different materials or they tend to be bigger.
For this reason the THz region is also sometimes called the quasi-optical or Gaussian region.
It is called the Gaussian region because beams in the THz region almost always are of Gaussian shape.
This, in turn, means that in the THz region the illumination conditions are non-homogeneous with an exponential decrease of illuminating signal strength towards the edges of the field-of-view.
Due to the large coherence length in this region of the electromagnetic spectrum, one can not simply make the illuminating beam larger and use only the center portion for imaging.
This will cause interference effects that appear as dominating artifacts in images (see Figure \ref{fig:artifacts}), which limit image quality and spatial resolution.
So, without calibration one either has to live with a limited field-of-view or one has to accept image artifacts and limited spatial resolution.
Therefore, calibration for non-homogeneous illumination is essential for a practical THz single-pixel imaging system.
We will discuss this issue and possible solutions in the next section.
Note that for a practical imaging system the calibration of the illuminating beam is very important and also leads to challenges related to the used masks.
This is also exemplified in Figure \ref{fig:convolution}, which shows a \SI{0.35}{\tera\hertz} single-pixel camera measurement of a non-metallic Siemenssstar test target reconstructed using a convolution approach with a non-negative least squares approach, i.e. we use convolutional masks leading to 2D circulant matrices.

In a modulated illumination setting of a THz single-pixel camera the quality and fidelity of the masks/illumination patterns determine the achievable spatial resolution, the signal-to-noise ratio, the achievable undersampling ratio and even more.
The physical process of implementing the masks is, therefore, very important and potentially introduces deviations into the masks already on the physical level.
These potential deviations can be simple blurring effects, the introduction of an offset or the reduction in the so-called modulation depth.
As introduced, for binary masks the modulation depth is essentially the difference between ones and zeros.
So, the ideal value is of course, 1 or 100\% but in a practical system the modulation depth can be several 10 percent below the ideal value.
Depending on the chosen reconstruction approach this will severely influence the image fidelity of reconstructed images (see Figure \ref{fig:convolution} for an example).
The example in Figure \ref{fig:convolution} shows what happens when the modulation depth is only 40\% and the masks are not optimized in an undersampling modality.
Due to the fact that the spatial resolution and the signal-to-noise ratio in the image are severely limited the image appears blurry and noisy.

All the aforementioned issues have to be considered on the software level in order to harness the potential power of a THz single-pixel camera and, therefore, often robustness is an important consideration when choosing the reconstruction method.
As mentioned the reconstruction approach used in this example was based on the idea of convolutional mask, which offers strong simplifications in the design of masks and the memory consumptions, moreover the reconstruction algorithms can be made more efficient.
Hence, such an approach has a lot of potential for THz single-pixel cameras, but there is still a lot of effort necessary in order to optimize the imaging process.

\section{Calibration Problems}
 
In the previous section we have seen the difficulties to calibrate the illumination beam directly, hence it seems necessary to perform a self-calibration approach during the reconstruction.
We hence look for a multiplicative decomposition of the image $\vimage$ into a smooth light field $\vlight$ and a normalized target structure $\vtarget$, i.e. 
\begin{equation}
\vimage = \vlight\odot\vtarget,
\end{equation} 
where $\odot$ denotes pointwise multiplication.
A simple approach would be to first reconstruct the image and then use standard decomposition ideas.
However, this approach seems suboptimal for undersampled measurements, since better a-priori knowledge on light and target is available than for their composition.

\subsection{Self-Calibration and Bilinear Inverse Problems}
A very recent trend in compressed sensing is to include the
calibration task into the recovery step as well.
This approach has
been termed also as {\em self-calibration}.
In the case of unknown
illumination this yields a {\em bilinear inverse problem}.
In
particular, for sparsity promoting regularizers this falls in the
category of biconvex compressed sensing \cite{Ling2015}.
In this work
the lifting approach has been adopted to transform the task to an
convex formulation, see here also the seminal work on Phaselift \cite{candes2013phaselift}.
Unfortunately, this approach does scale for
imaging problems.

Here we are confronted with the generic problem:
\begin{equation}
 \min_{\vtarget\in [0,1]^n,\vlight\geq 0}
 \|\meas-\Meas(\vlight\odot\vtarget)\|^2_{\ellp{2}}+ \lambda\cdot r(\vtarget,\vlight)
\end{equation}
with an appropriate regularization functional $r$ for both structures,
the illumination $\vlight$ and the target $\vtarget$.
This problem is
a particular case of a {\em bilinear inverse problem} and linked to
{\em compressive blind deconvolution}.
Indeed, let us formalize this
by using the 2D fast Fourier transform (FFT):
\begin{equation}
 \vlight\odot\vtarget=\frac{1}{\sqrt{n}}\Fourier^{-1}[(\Fourier\vlight)\ast(\Fourier\vtarget)]
 =\frac{1}{\sqrt{n}}\Fourier^{-1}[\vlightLog\ast\vtargetLog]
\end{equation}
where $\Fourier\colon\mathbb{R}^n\rightarrow\mathbb{C}^n$ defines the Fourier transform operator and $\vlightLog,\vtargetLog\in\mathbb{C}^n$.
and therefore the observation
$\meas=\frac{1}{\sqrt{n}}\Meas\Fourier^{-1}(\vlightLog\ast\vtargetLog)$
is a compressed circular 2D convolution of $\vlightLog$ and $\vtargetLog$.
Obviously, without further assumptions this type of inverse problem
cannot be solved uniquely.

Let us discuss the case $\Meas=\Id$ that usually corresponds to an
image scanning approach (and not the single pixel setup).
Here, the measurement image is $\vimage$ and
the goal is to factorize it into light $\vlight$ and a normalized target $\vtarget$
using the program:
\begin{equation}
 \min_{\vtarget\in [0,1]^n,\vlight\geq 0}
 \|\vimage-\vlight\odot\vtarget\|^2_{\ellp{2}}+ \lambda\cdot r(\vtarget,\vlight).
\end{equation}
Note that this formulation is a non-convex problem since it is not
jointly convex in $(\vlight,\vtarget)$.
To make this problem more well-posed from a compressed viewpoint further assumptions on $\vlight$ and $\vtarget$ are necessary.
In the case of random subspace assumptions for either $\vlight$ or $\vtarget$ and some additional incoherence conditions, recovery guarantees for the convexified (lifted) formulation of blind deconvolution have been established in \cite{Ahmed:2012}.
This framework even extends to the more involved problem of blindly demixing convolutional mixtures, called as also blind demixing and deconvolution, and recently almost-optimal sampling rates could established here under similar assumptions
\cite{Jung:TIT:2017}.
However, the analysis of the original non-convex problem (without lifting) is often difficult, requires a priori assumptions and sufficient randomization and the performance of iterative decent algorithms often depends on a good initialization.
First results in this direction appeared here for example in \cite{Li2016}.

A further way to convexify a related problem is based
on a formulation for strictly positive gains
\begin{equation}
 \min_{\vtarget\in [0,1]^n,\vc{g}>\epsilon}
 \|\vimage-\vc{g}\odot\Meas\vtarget\|^2_{\ellp{2}}+ \lambda\cdot r(\vtarget,\vc{g}).
\end{equation}
In this problem the image $\vtarget$ is unknown and each measurement
has an unknown strictly-positive gain $\vc{g}>\epsilon$.
This problem
occurs exactly in our setting when $\Meas=\Id$ and then $\vc{g}=\vlight$.
Thus, under the additional prerequiste of $\vlight\geq\epsilon>0$ one could
write a constraint $\vimage=\vlight\odot\vtarget$ also as
$\vlight^{-1}\odot\vimage=\vtarget$ which is jointly convex in
$(\vlight^{-1},\vtarget)$ \cite{gribonval2012blind} \cite{ling2018self}.
This approach
is also interesting if $\vtarget=\Dict\vc{\vcoeff}$ itself is a
compressed representation.
Unfortunately, this approach does not apply to the single pixel imaging
setup due to matrix $\Meas$, which cannot be interchanged with the division by $\vlight$.

\subsection{Single Pixel and Retinex Denoising}

A related classical problem in imaging is the so-called Retinex approach.
The corresponding problem was first investigated by Land (cf. \cite{land1964retinex} \cite{Land:retinex}) in the context of human visual studies and provide a first entry in the developing Retinex theory.
The human eye is able to compensate lack of illumination, hence it is able to filter for structural information in the field of view.
From a more general perspective this translates to the question of how to separate an image into a structural and an illumination part.
Here we focus on the approach of Kimmel et al. \cite{Kimmel:2003} and further developments of this approach.
By defining $\vimageLog \defgr \log(\vimage)$, $\vtargetLog \defgr \log(\vtarget)$ and $\vlightLog \defgr \log(\vlight)$ we move the image into the logarithmic domain.
Again in the case $\Meas=\Id$ this allows for the usage of a convex variational model as proposed in \cite{Kimmel:2003}.
The basic assumptions are as follows:
\begin{itemize}
	
	\item {\em spatial smoothness of illumination},
	\item {\em $\vlightLog$ is greater than $\vimageLog$}: Since $\vtarget\in [0,1]^n$ and $\vlight \geq \vtarget\odot\vlight$ it is $\vlightLog\geq\vimageLog$ due to the monotone nature of the logarithmic map,
	\item {\em non trivial illumination}: resemblance to the original image, hence $l-s$ small,
	\item {\em soft smoothing on structure reconstruction} and
	\item {\em smooth boundary condition}.
\end{itemize}
These result in the variation approach proposed by Kimmel et al.:
\begin{equation}
	\min_{\vlightLog \geq \vimageLog} \frac{1}{2} \|\vlightLog-\vimageLog\|^2_{\ellp{2}} + \frac{\alpha}{2} \|\nabla(\vlightLog-\vimageLog)\|^2_{\ellp{2}} + \frac{\beta}{2}\|\nabla \vlightLog\|^2_{\ellp{2}}.\label{eq:kimmel}
\end{equation}
However, this applies a smoothing on $l$ and $l-s$, which is basically $r$.
Then this boils down to the data fitting term and two smoothing regularization operators that are balanced using the respective parameters $\alpha$ and $\beta$.
In a later paper by Ng and Wang \cite{ng_total_2011}, the usage of TV regularization for the structural image was introduced.
The proposed model is
\begin{equation}
	\min_{\vlightLog\geq \vimageLog, \vtargetLog\leq 0} \frac{1}{2}\|\vlightLog-\vimageLog+\vtargetLog\|^2_{\ellp{2}} + \alpha \|\nabla \vtargetLog\|_{\ellp{1}} + \frac{\beta}{2}\|\nabla \vlightLog\|^2_{\ellp{2}}.\label{eq:ngandwang}
\end{equation}
This equation makes sense from a Retinex point of view.
One has given an image $\vimage$ or $\log(\vimage)=\vimageLog$ respectively, and wants to separate reflection $\log(\vtarget)=\vtargetLog$ and $\log(\vlight)=\vlightLog$.

\begin{figure}[t]
	\centering
	\begin{tabular}{B{0.23\textwidth}B{0.23\textwidth}B{0.04\textwidth}B{0.23\textwidth}B{0.23\textwidth}}
		\multicolumn{2}{c}{L2 based (see \eqref{eq:kimmel})}&&\multicolumn{2}{c}{TV based (see \eqref{eq:ngandwang})}\\
		\includegraphics[width=\linewidth]{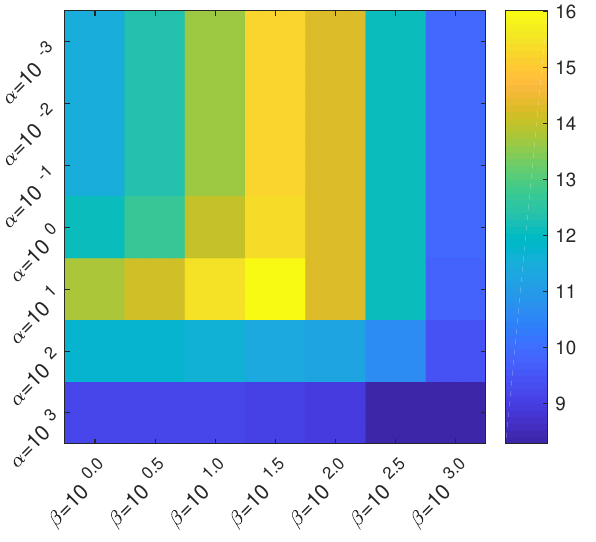}&
		\includegraphics[width=\linewidth]{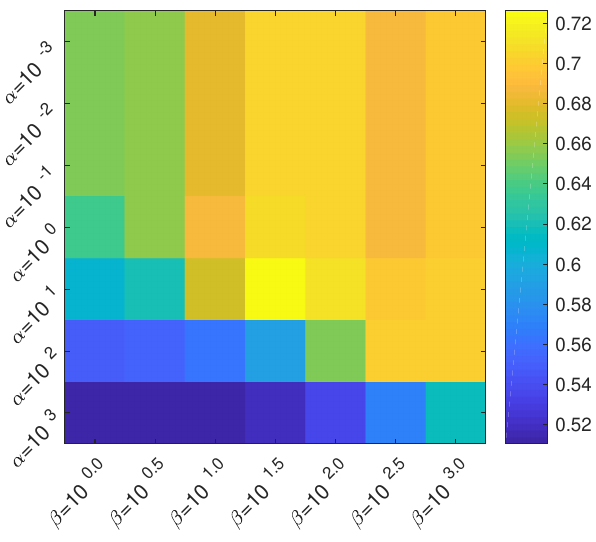}&&
		\includegraphics[width=\linewidth]{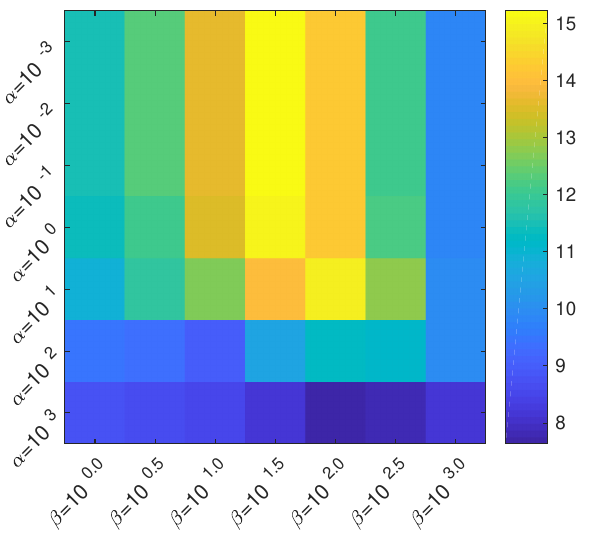}&
		\includegraphics[width=\linewidth]{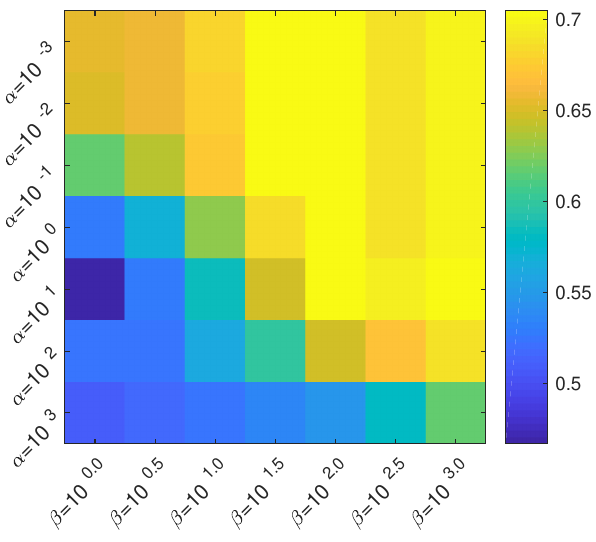}\\
		PSNR&SSIM&&PSNR&SSIM\\
	\end{tabular}
	\caption{Overview of PSNR value and SSIM index for systematic parameter test for the synthetic x phantom dataset.
	Map of regularization parameters in respect to given pair of regularization parameters $\alpha$ and
$\beta$: $\alpha$ in $Y$ axes with values $10^{-3}$ (top), $10^{-2}$, ..., $10^2$, $10^3$ (bottom) and $\beta$ in $X$ axis
with values $10^0$ (left), $10^{0.5}$, ..., $10^{2.5}$, $10^3$ (right).}
	\label{im:comparePrametertest}
\end{figure}

\begin{figure}[t]
	\centering
	\includegraphics[width=\linewidth]{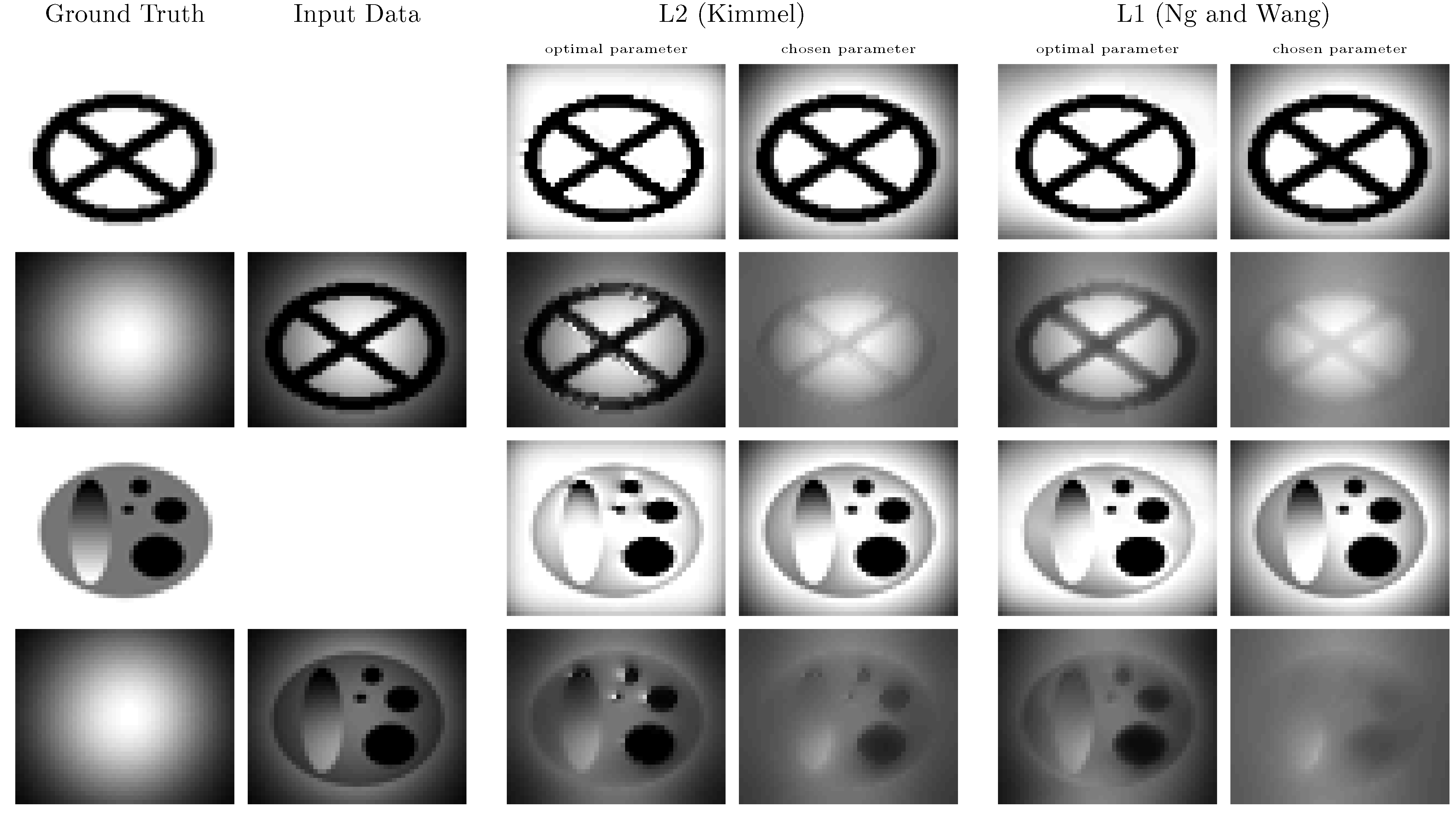}
	\caption{Comparison of the L2 based Retinex model by Kimmel and the total variation based Retinex model by Ng and Wang for both datasets.}
	\label{im:retinexComparison}
\end{figure}

A comparison of both Retinex models are shown in Figure \ref{im:retinexComparison}, with an interesting parameter dependent behavior.
The regularization and weighting parameters $\alpha$ and $\beta$ have been determined using a parameter test, shown in Figure \ref{im:comparePrametertest}.
We see that choosing parameters optimizing the SSIM measure for the target respectively its logarithm $\vtargetLog$ ($\alpha=10^0$, $\beta=10^2$ in the Kimmel model respectively $\alpha=10^1$, $\beta=10^2$ for the Ng and Wang model) yields a strange result with respect to the illumination, which still has quite some of the target structure included, in addition to the full beam.
This can be balanced by different parameter choices ($\alpha=10^0$, $\beta=10^3$ in the Kimmel model respectively $\alpha=10^1$, $\beta=10^3$ for the Ng and Wang model) displayed in the lower line respectively, which eliminates most structure from the reconstructed illumination, but also leaves some beam effects in the target.


\begin{figure}[t]
	\centering
	\includegraphics[width=\linewidth]{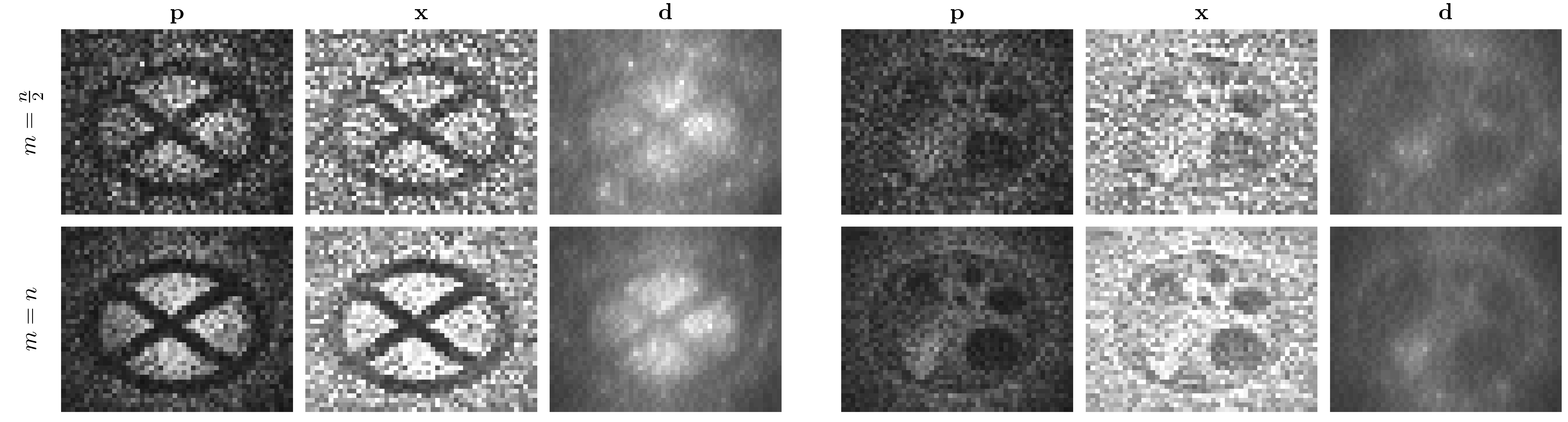}
	\caption{Reconstruction results for an alternating reconstruction and Retinex scheme.}
	\label{im:compareReconRetinex}
\end{figure}

However in the framework of single pixel camera approaches, one does not have the fully recovered image at hand.
Instead we only have measured data, hence the data fidelity term in \eqref{eq:ngandwang} is replaced by a reconstruction based data fidelity term based on the forward operator $A$ and measured data $y$:
\begin{equation}
 \min_{\vlightLog\geq \vimageLog, \vtargetLog\leq 0}
 \|\meas-\Meas(e^{\vlightLog+\vtargetLog})\|^2_{\ellp{2}}+
 \alpha\|\nabla \vtargetLog\|_{\ellp{1}}+\frac{\beta}{2}\|\nabla\vlightLog\|^2_{\ellp{2}}.
\end{equation}
A simple approach is to apply a two-step idea: In the first step one can compute standard reconstruction of $\vimage$ e.g. by \eqref{eq:TVmodel} and subsequently apply the Retinex model with $\vc{s}=\log \vimage$.
Since in this case the first reconstruction step does not use prior knowledge about the structure of illumination, it is to be expected that the results are worse compared to a joint reconstruction when the number of measurements decreases.

\begin{figure}[t]
\centering
\includegraphics[width=\linewidth]{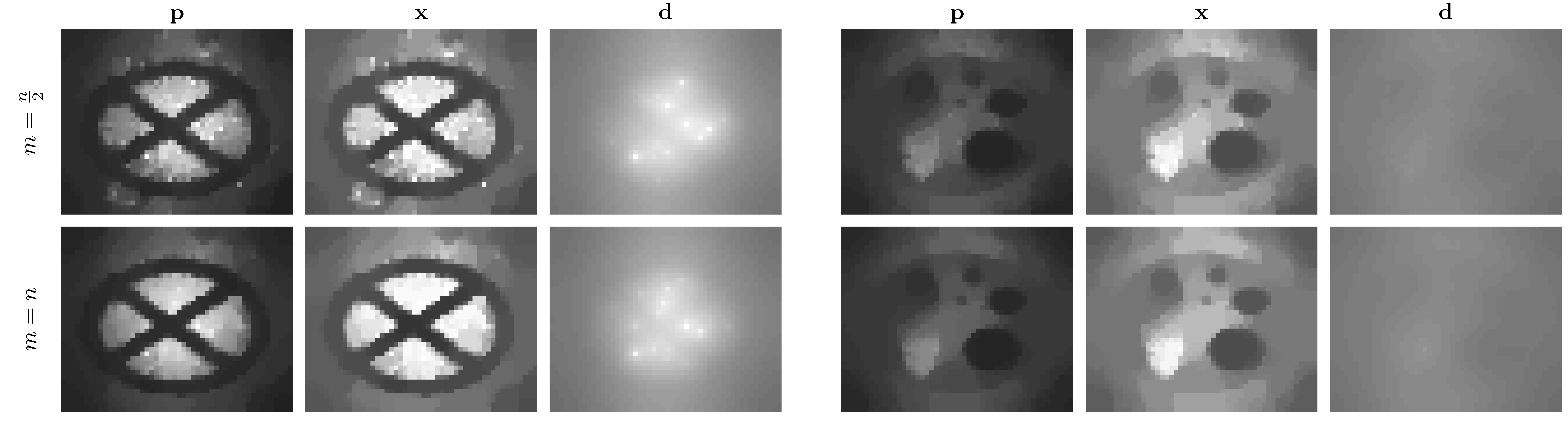}
\caption{Reconstruction results for Forward-Backward Splitting including the L1 based model proposed by Ng and Wang.}
\label{im:compareReconRetinexTV}
\end{figure}
In order to construct a computational approach for the calibration problem we still formulate the problem in the logarithmic variables
\begin{equation}
\min_{\vc{r}\leq 0,\vc{l}\in\mathbb{R}^n} \frac{1}2 \lVert\meas-\Meas e^{\vc{r}+\vc{l}} \rVert^2_{\ellp{2}} + \lambda \cdot \tilde{r}(\vc{r},\vc{l})
\end{equation}
and derive a forward-backward splitting algorithm in these variables: The forward step is simply given by
\begin{align*}
 \vc{r}^{k+1/2} &= \vc{r}^{k} + \tau e^{\vc{r}^k+\vc{l}^k} \odot \Meas^T (\meas-\Meas e^{\vc{r}^k+\vc{l}^k}) \\
 \vc{l}^{k+1/2} &= \vc{l}^{k} + \tau e^{\vc{r}^k+\vc{l}^k} \odot \Meas^T (\meas-\Meas e^{\vc{r}^k+\vc{l}^k}) 
\end{align*}
and the backward step then computes $\vc{r}^{k+1}$ as a minimizer of
$$ \frac{1}{2\tau} \Vert \vc{r} - \vc{r}^{k+1/2} \Vert_{\ell_2}^2 + \lambda_r \Vert \nabla r \Vert_{\ell_1}$$
and $\vc{l}^{k+1}$ as a minimizer of
$$ \frac{1}{2\tau} \Vert \vc{l} - \vc{l}^{k+1/2} \Vert_{\ell_2}^2 + \lambda_l \Vert \nabla l \Vert_{\ell_2}^2. $$

Note that by adding the two equations in the forward step we can directly formulate it as an update for $\vc{s}$ as 
$$ \vc{s}^{k+1/2} = \vc{s}^{k} + \tau e^{\vc{s}^k} \odot \Meas^T (\meas-\Meas e^{\vc{s}^k }) .$$
This induces a simple idea for the two-step approach: we can iterate $\vc{s}$ with the above scheme and directly apply the 
above Retinex model to data $\vc{s}^{k+1/2} $.
From the resulting minimizers $\vc{r}^{k+1}$ and $\vc{l}^{k+1}$ we can compute 
$ \vc{s}^{k+1} = \vc{r}^{k+1} + \vc{l}^{k+1}$.
The resulting reconstructions are shown in Figure \ref{im:compareReconRetinex}, which are clearly suboptimal since too much of the random structure from the matrix $\Meas$ is propagated into $\vc{s}^{k+1/2} $, which is not reduced enough in the Retinex step.
The results of the forward-backward splitting approach are shown in Figure \ref{im:compareReconRetinexTV}, which are clearly improving the separation of illumination beam and target and yield robustness with respect to undersampling, although still not providing perfect smoothing of the images.
This might be expected from improved regularization models for the decomposition to be investigated in the future.

\section{Phase Retrieval in Single Pixel Cameras}

In the following we discuss another aspect of reconstruction in single
pixel cameras, namely phase reconstruction problems caused by
diffraction effects.
This problem is a particular instance of the
difficult phase retrieval problem.
Compared to blind deconvolution, as
a (non-convex) bilinear inverse problem, phase retrieval is the
corresponding quadratic (non-convex) case and therefore most analytical results
here are based on lifting the problem to a convex formulation
\cite{candes2013phaselift}.
For some overview, further references and
due to limited space we refer here exemplary to the overview article
\cite{Jaganathan:phaseretrieval2016}.
Interestingly also, first works already
appeared where phase retrieval and blind deconvolution are combined
using lifting \cite{Ahmed:nips18}.
However, already for imaging
problems at moderate resolution these approaches usually not scale and
require a priori random subspace assumptions which are difficult to
fulfill in practice.
In the following we will discuss how to setup the
phase reconstruction problem in the diffraction case.

The propagation of light waves can - after some approximation - be represented by the Fresnel diffraction integral, we refer to \cite{Nickel} and references cited therein for a detailed treatment of the discretization problem.
In the single pixel setup (compare Figure \ref{diffractionfigure}) we consider that the diffraction has to be taken into account in two ways, namely between the object and mask plane and the mask and detector plane.
This means that the measurement matrix is further changed by the introduction of the masks.
In addition to this only the magnitude of the combined complex signals is obtained.
Let ${\bf D}_{om}$ be a matrix discretizing the diffraction integral from object to mask plane and ${\bf D}_{md}$ be a matrix discretizing the diffraction from mask to detector plane.
Then the complex signal arriving at the detector plane is
$$ {\bf z}_i = {\bf 1}^T {\bf D}_{md} \text{ diag}(\Meas_i) {\bf D}_{om} \vimage, $$
where $\vimage$ is the complex image, $\Meas_i$ denotes the $i$-th row of $\Meas$, and ${\bf 1}$ a vector filled with all entries equal to one.
Note that in the absence of diffraction, i.e. ${\bf D}_{md} = {\bf D}_{om} = {\bf Id}$, this reduces to the standard single pixel camera approach discussed above.
The measured intensity is then the absolute value of ${\bf z}_i$, or rather its square, i.e. 
$$ \meas_i = \vert {\bf z}_i \vert^2 = \vert {\bf 1}^T {\bf D}_{md} \text{ diag}(\Meas_i) {\bf D}_{om} \vimage \vert^2. $$

Defining a matrix ${\bf B}$ with rows 
$$ {\bf B}_i = {\bf 1}^T {\bf D}_{md} \text{ diag}(\Meas_i) {\bf D}_{om}, $$
the phase reconstruction problem can be written in compact notation as
\begin{equation} \label{phaseretrieval}
\meas = \vert {\bf B} \vimage \vert^2.
\end{equation} 
It is apparent that \eqref{phaseretrieval} is a nonlinear problem compared to the linear phase- and diffractionless problem, where both the matrix $\Meas$ and the image $\vimage$ can be modeled to be real non-negative (hence the absolute value does not lead to a loss of information).
For this reason the iterative solution of \eqref{phaseretrieval} respectively least-squares versions thereof is a problem of central importance, which we will discuss in the next section.

\subsection{Algorithms for Phase Retrieval}

Phase retrieval is a classical problem in signal processing, which has been revived by compressed sensing approaches in the last decade, 
cf. \cite{fienup1982phase} \cite{bauschke2002phase} \cite{marchesini2007invited} \cite{jaganathan2015phase} for an overview of classical and recent methods.
A famous early algorithm is the Gerchberg-Saxton (GS) algorithm (cf. \cite{gerchberg1972practical}).
However a version of the GS algorithm has been developed by Fienup (cf. \cite{fienup1982phase}) that works with amplitude measurements only and is using a multiplicative update scheme of the form
$$ \meas_{k+1} = \frac{\vert \meas \vert}{\vert {\bf B} \vimage_k \vert} ~{\bf B} \vimage_k , \quad \vimage_{k+1} = {\bf B}^+ \meas_{k+1},$$
where ${\bf B}^+$ defines the commonly used pseudo inverse of ${\bf B}$ (cf. \cite{engl_regularization_2000}).
In the case of additional constraints on $\vimage$, those are applied to modify $\meas_{k+1}$ in an additional projection step.
Note that the original Gerchberg-Saxton algorithm is formulated for Fourier measurements only, where ${\bf B}$ is invertible by its adjoint.
The algorithm can be used in particular to compute real images, where 
$$ \vimage_{k+1} = \text{Re }({\bf B}^+ \meas_{k+1}) .$$

Another approach that can be found at several instances in literature (cf. \cite{blaschke1997regularization} \cite{shechtman2014gespar} \cite{gao2017phaseless}) is based on the application of (variants of) Gauss-Newton methods to the 
least squares problem of minimizing
\begin{equation}
L(\vimage) = \Vert \meas - \vert {\bf B} \vimage \vert^2 \Vert^2.
\end{equation} 
This amounts to linearizing the residuals around the last iterate and to obtain $\vimage_{k+1}$ as the minimizer of
\begin{equation}
\sum_{i=1}^m \vert \meas_i - \vert {\bf B}_i \cdot \vimage_k \vert^2 + 2(\text{Re}({\bf B_i}) \text{Re}({\bf B_i}) \cdot \vimage_k + \text{Im}({\bf B_i}) \text{Im}({\bf B}_i) \cdot \vimage_k)\cdot (\vimage - \vimage_k )\vert^2.
\end{equation}
Several variants are used to stabilize the Gauss-Newton iteration and to account for the ill-conditioning of the Jacobian matrix 
$$ {\bf J}_k = 2\left((\text{Re}({\bf B_i}) \text{Re}({\bf B_i}) \cdot \vimage_k + \text{Im}({\bf B_i}) \text{Im}({\bf B_i}) \cdot \vimage_k) \right)_{i=1,\ldots,m}. $$
A popular one, which we will also use in our numerical tests, is the Levenberg-Marquardt method, which computes $\vimage_{k+1}$ as the minimizer of 
\begin{equation}
 \Vert \meas - \vert {\bf B} \vimage_k \vert^2 + {\bf J}_k (\vimage - \vimage_k )\Vert^2 + \alpha_k \Vert \vimage - \vimage_k \Vert^2,
\end{equation}
with $\alpha_k$ a decreasing sequence of positive parameters.

A recently popular approach to solve the nonconvex least-squares problem is the Wirtinger flow (cf. \cite{candes2015phase}).
Its main ingredient is just gradient descent on the least squares functional $L$, i.e., 
\begin{align}
\vimage_{k+1} &= \vimage_k - \frac{\mu_{k+1}}{2 m\Vert \vimage_0 \Vert^2 } \nabla L(\vimage_k) \\ &= 
\vimage_k + \frac{\mu_{k+1}}{ m\Vert \vimage_0 \Vert^2 } \sum_{i=1}^m (\meas_i - \vert {\bf B}_i \cdot \vimage_k \vert^2) (\text{Re}({\bf B}_i {\bf B}_i^T)+ \text{Im}({\bf B}_i {\bf B}_i^T)). \nonumber
\end{align} 
For the choice of the step size an increasing strategy like $\mu_k \sim 1- e^{-k/k_0}$ is proposed.
Obviously the gradient descent as well as the Levenberg-Marquardt method above rely on appropriate initializations in order to converge to a global minimum.
For this sake a spectral initialization has been proposed, which chooses $\image_0$ as the first eigenvector of the positive semidefinite matrix 
$${\bf M}_0 = \sum_{i=1}^m \meas_i (\text{Re}({\bf B}_i {\bf B}_i^T)+ \text{Im}({\bf B}_i {\bf B}_i^T), $$ 
corresponding to the data sensitivity in the least squares functional.
For practical purposes the first eigenvector can be approximated well with the power method.

A very recent approach is the truncated amplitude flow (cf. \cite{wang2018solving}), whose building block is gradient descent on the alternative least squares functional
\begin{equation}
\tilde L(\vimage) = \sum_{i=1}^m \vert \sqrt{\meas_i} - \vert {\bf B}_i^T \vimage \vert \vert^2,
\end{equation} 
which is however not differentiable for ${\bf B}_i^T \vimage = 0$.
In the case the derivative exists we find
$$ \tilde L(\vimage) = - \sum_{i=1}^m( \sqrt{\meas_i} - \vert {\bf B}_i^T \vimage \vert ) \frac{1}{\vert {\bf B}_i^T \vimage \vert}
{\bf B}_i {\bf B}_i^T \vimage. $$ 
The truncated amplitude flow now only selects a part of the gradient in order to avoid the use of small $\vert {\bf B}_i^T \vimage \vert$ (compared to $\sqrt{\meas_i}$), i.e. 
\begin{equation}
\vimage_{k+1} = \vimage_k + \mu_k \sum_{i\in I_k}\left( \frac{\sqrt{\meas_i}}{\vert {\bf B}_i^T \vimage_k \vert} - 1 \right) 
{\bf B}_i {\bf B}_i^T \vimage_k,
\end{equation}
with index set
$$ I_k = \left\{ i \in \{1,\ldots,m\}~|~\frac{\sqrt{\meas_i}}{\vert {\bf B}_i^T \vimage_k \vert} \leq 1+ \gamma \right\} $$
with some $\gamma > 0$.
For the truncated amplitude flow another initialization strategy has been proposed in \cite{wang2018solving}, which tries to minimize the angles of $\vimage_0$ to the measurement vectors ${\bf B}_i$ for a subset of indices $I_0$.
This is equivalent to computing the eigenvector for the smallest eigenvalue of the matrix 
$${\bf \tilde M}_0 = \sum_{i\in I_0} \frac{1}{\Vert {\bf B}_i \Vert^2} {\bf B}_i {\bf B}_i^T. $$
Since computing such an eigenvector is a problem of potentially high effort, it is approximated by computing the largest eigenvalue of a matrix built of the remaining rows of ${\bf B}$ (cf. \cite{wang2018solving} for further details).

In order to introduce some regularization into any of the flows above, we can employ a forward-backward splitting approach.
E.g. we can produce the Wirtinger flow to produce the forward estimate and use an additional backward proximal step
\begin{align}
\vimage_{k+1/2} &= \vimage_k - \frac{\mu_{k+1}}{2 m\Vert \vimage_0 \Vert^2 } \nabla L(\vimage_k) \\
\vimage_{k+1} &= \text{prox}_{\lambda_{k+1} R}(\vimage_{k+1/2})
\end{align}
with the regularization functional $R$ and a parameter $\lambda_{k+1}$ chosen appropriately in dependence of $\mu_{k+1}$.
The proximal map 
$\vimage = \text{prox}_{\lambda_{k+1} R}(\vimage_{k+1/2})$ is given as the unique minimizer of 
$$ \frac{1}2 \Vert \vimage - \vimage_{k+1/2} \Vert^2 + \lambda_{k+1} R(\vimage). $$

Finally we mention that in addition to the nonconvex minimization approaches there has been a celebrated development towards convexifying the problem, the so-called PhaseLift approach, which can be shown to yield an exact convex relaxation under appropriate conditions (cf. \cite{candes2013phaselift}).
The key idea is to embed the problem into the space of $n \times n$ matrices and to find $\vimage \vimage^T$ as a low rank solution of a semidefinite optimization problem.
The squared absolute value is rewritten as 
$$\vert {\bf B}_i \vimage\vert^2 = {\bf B}_i^* {\bf P} {\bf B}_i$$
with the positive semidefinite rank-one matrix ${\bf P} = \vimage \vimage^*$.
The system of quadratic equations $\vert {\bf B}_i \vimage\vert^2=y_i$ then is rewritten as the problem of finding the semidefinite matrix of lowest rank that solves the linear matrix equations ${\bf B}_i^* {\bf P} {\bf B}_i = y_i$.
Under appropriate condition this problem can be exactly relaxed to minimizing the nuclear norm of ${\bf P}$, a convex functional, subject to linear equation and positive semidefiniteness.
 The complexity of solving the semidefinite problem in dimension $n^2$ makes this aprroach prohibitive for applications with large $n$ however, hence we will not use it in our computational study in the next section.

\subsection{Results}

\begin{figure}[t]
	\centering
	\includegraphics[width=.8\textwidth]{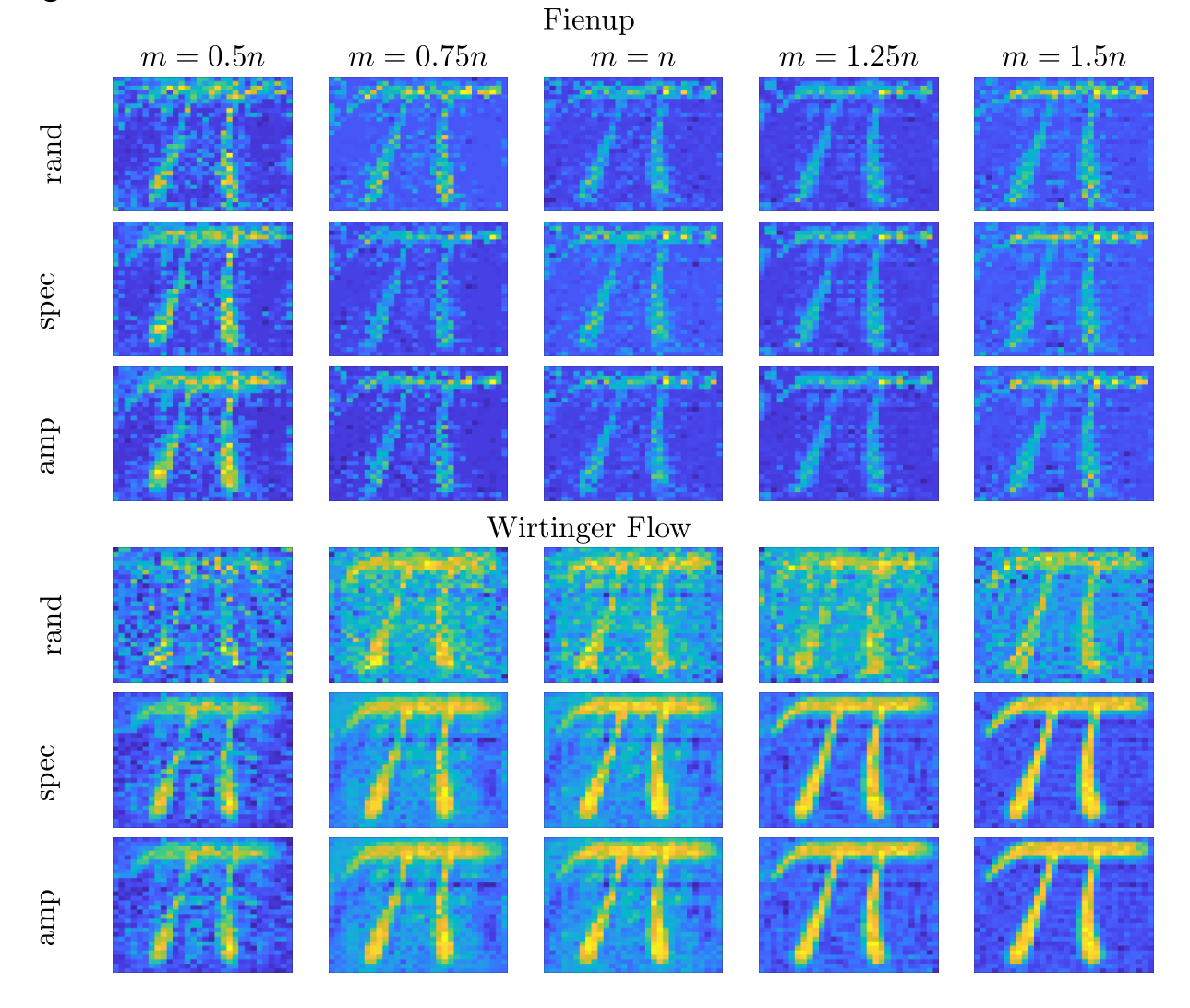}
	\caption{Reconstructions of $\pi$-image with the Fienup-variant of the Gerchberg-Saxton algorithm and the Wirtinger flow (from \cite{Nickel}). \label{reconfigure1}}
\end{figure}

\begin{figure}[t]
	\centering
	\includegraphics[width=.8\textwidth]{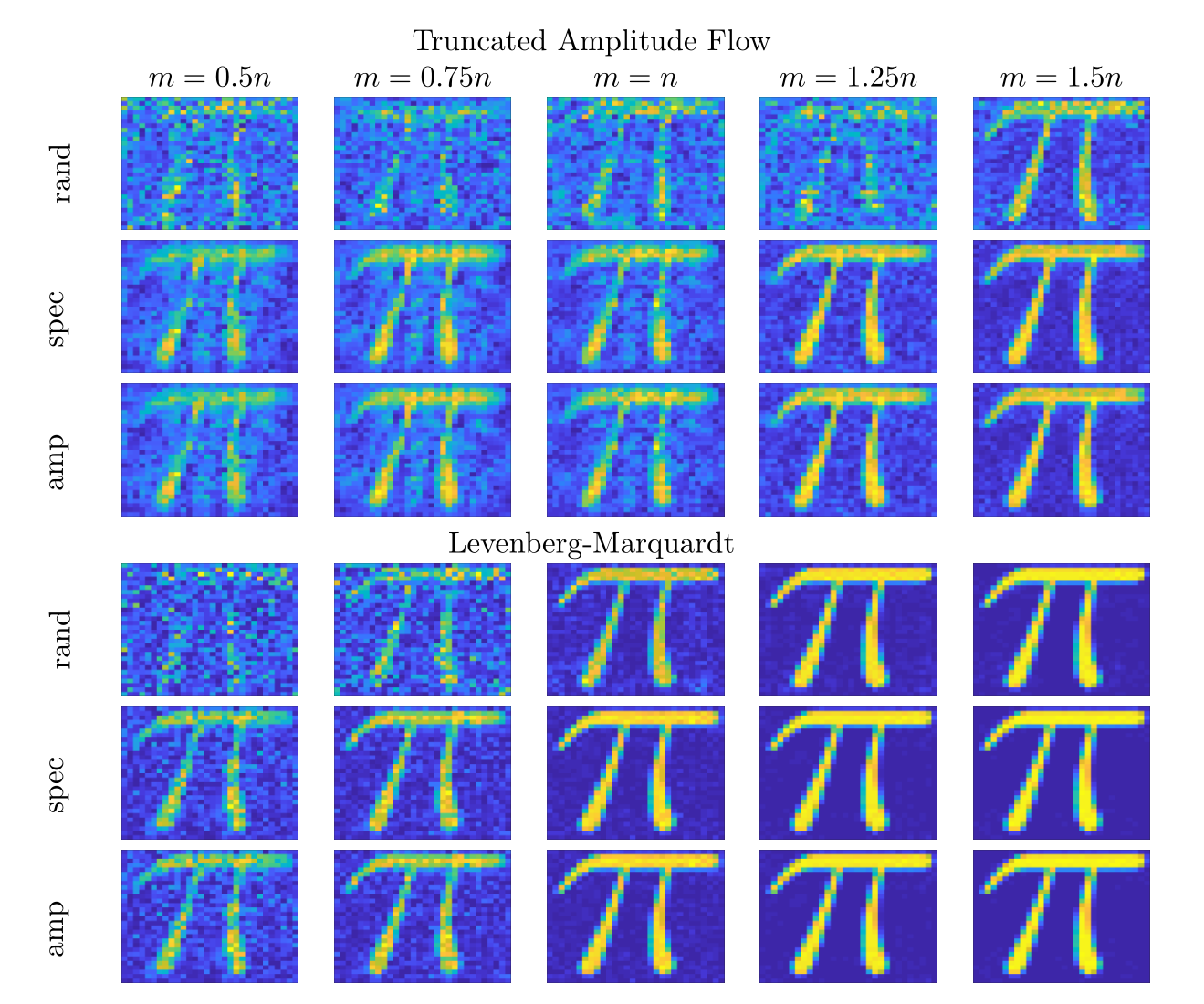}
	\caption{Reconstructions of $\pi$-image with the truncated amplitude flow and the Levenberg-Marquardt method (from \cite{Nickel}). \label{reconfigure2}}
\end{figure}

\begin{figure}[t]
	\centering
	\includegraphics[width=0.75\textwidth]{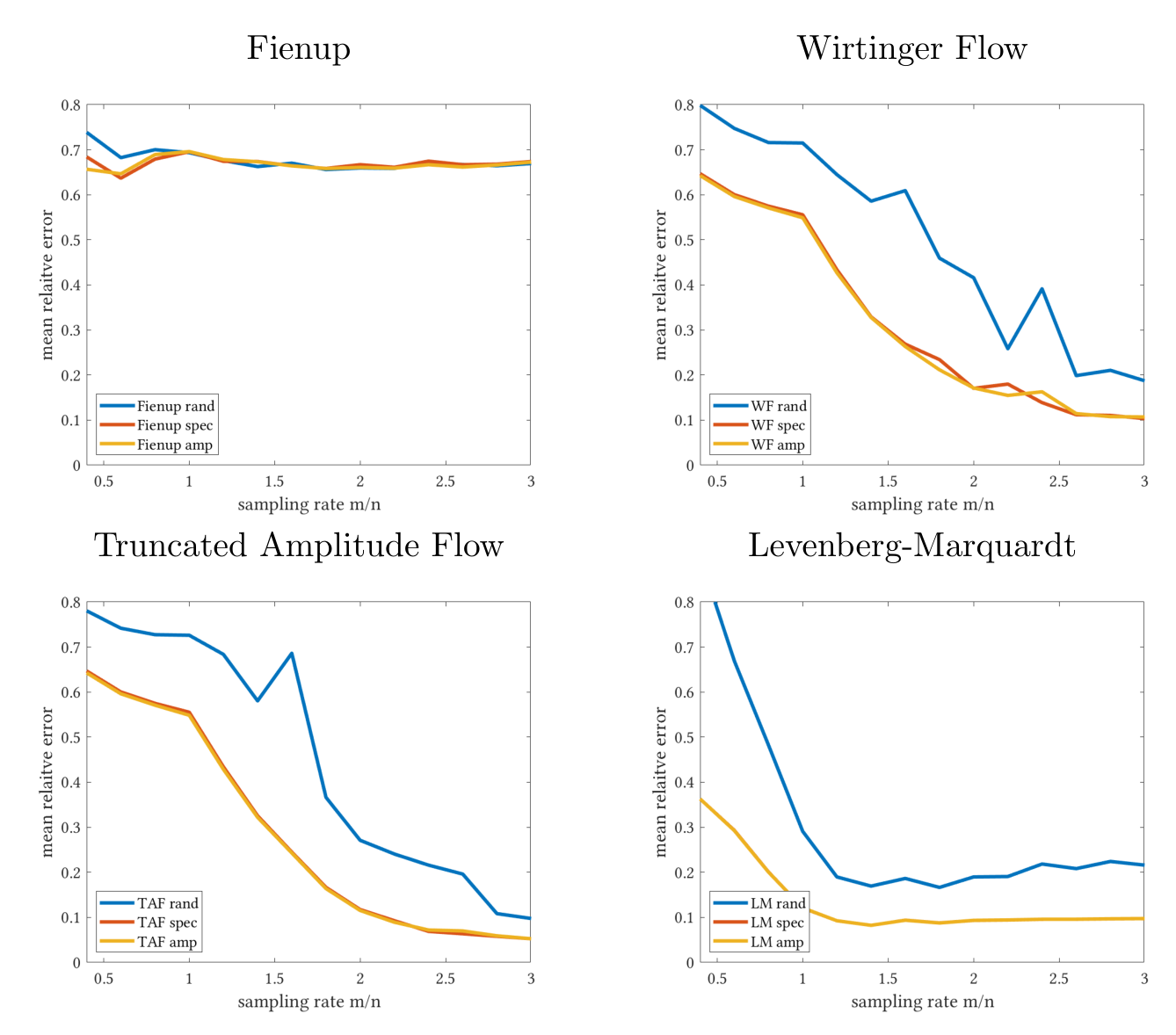}
	\caption{Plot of relative mean square error to ground truth vs. sampling rate for reconstructions of $\pi$-image (from \cite{Nickel}). \label{reconfigure3}}
\end{figure}

In the following we present some results obtained from the algorithms discussed above when applied to our single pixel setup and compare their performance.
We start by using the algorithms from the previous section for reconstructing a real-valued image showing the figure $\pi$ ($n=900$) without noise in the measurements.
For this sake we consider different types of undersampling and three different initializations: a random initial value (rand), the spectral initialization proposed originally for the Wirtinger flow (spec) and the orthogonality promoting initialization originally proposed for the truncated amplitude flow (amp).
The results for the four different methods are shown in Figures \ref{reconfigure1} and \ref{reconfigure2}.
The results clearly demonstrate the dependence of solutions on the initial value, in particular we see the improvement of the two new initialization strategies with respect to the random one (except in the case of the algorithm by Fienup, which however yields inferior results in all cases).
We observe that the Levenberg-Marquardt method performs particularly well for larger sampling rates and yields an almost perfect reconstruction, but also produces suitable reconstructions for stronger undersampling.
A more quantitative evaluation is given in Figure \ref{reconfigure3}, which provides plots of the relative mean-square error vs. the sampling rate, surprisingly the Levenberg-Marquardt method outperforms the other schemes for most rates.

\begin{figure}
	\centering
\includegraphics[width=.9\textwidth]{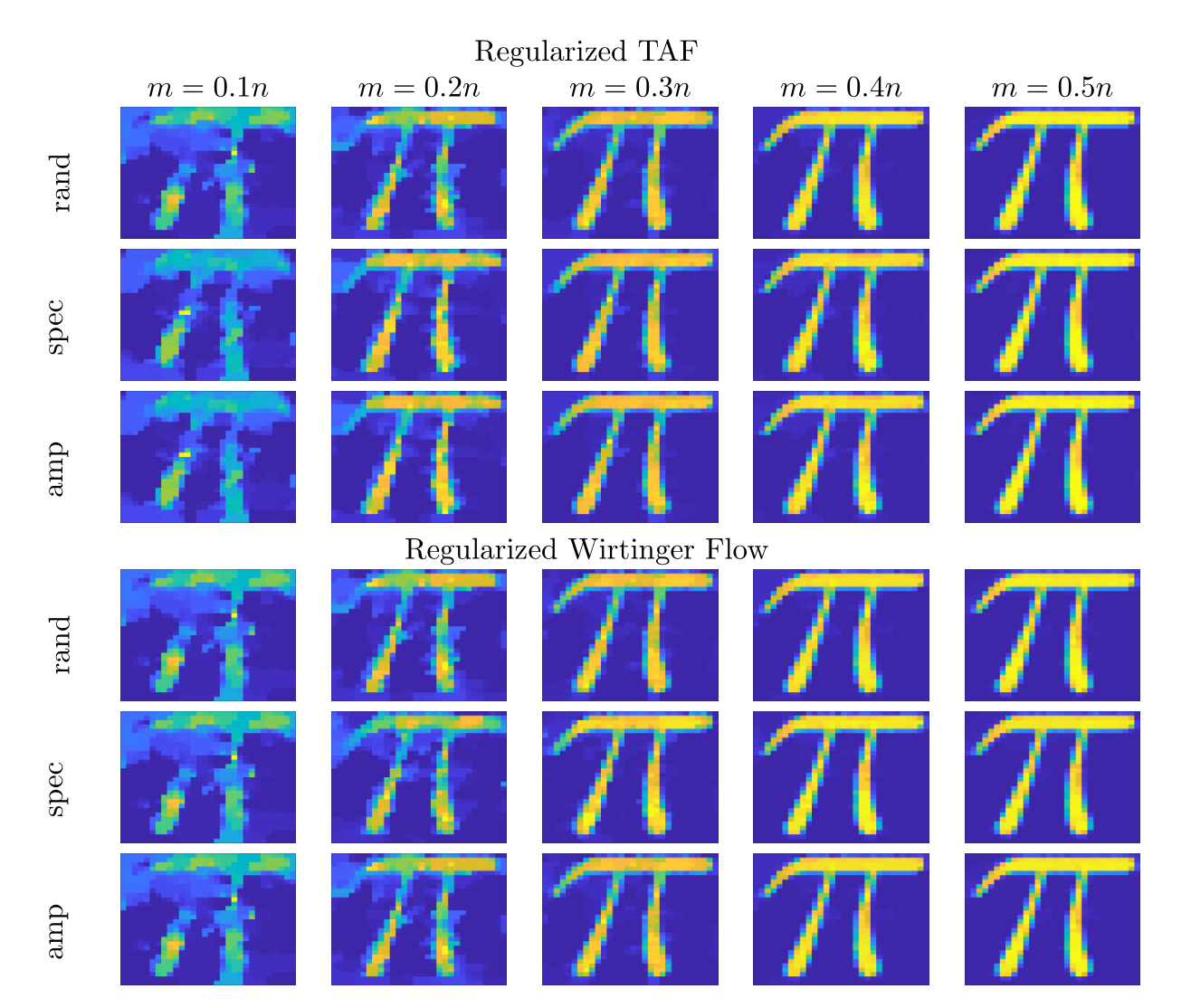}
\caption{Reconstructions of $\pi$-image with the regularized version of the truncated amplitude flow and the Wirtinger flow (from \cite{Nickel}). \label{regularizedFlows}}
\end{figure}

\begin{figure}
	\centering
\includegraphics[width=.9\textwidth]{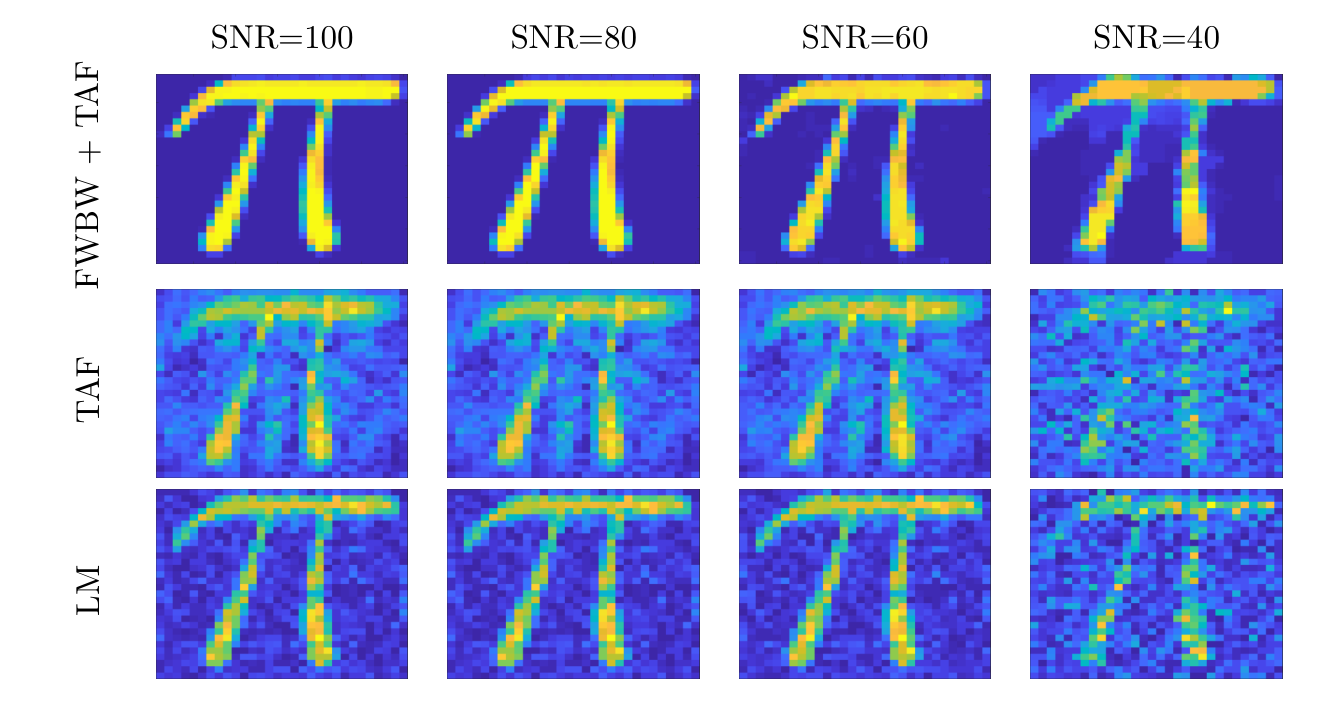}
\caption{Reconstructions of $\pi$-image from noisy data with different signal-to-noise ratios (SNR) (from \cite{Nickel}). \label{reconnoise}}
\end{figure}

The positive effect of regularization and its necessity for very strong undersampling is demonstrated in Figure \ref{regularizedFlows}, where we display the reconstruction results obtained with the forward-backward splitting strategy and the total variation as regularization functional.
We see that both approaches yield similar results and in particular allow to proceed to very strong undersampling with good quality results.
The effect of noise in the data, here for a fixed sampling ratio $\frac{m}n=0.7$ is demonstrated in Figure \ref{reconnoise} for different signal-to-noise ratios.
Again, not surprisingly, the regularized version of the flows yields stable reconstructions of good quality even for data of lower quality.

\begin{figure}[t]
	\centering
\includegraphics[width=0.75\textwidth]{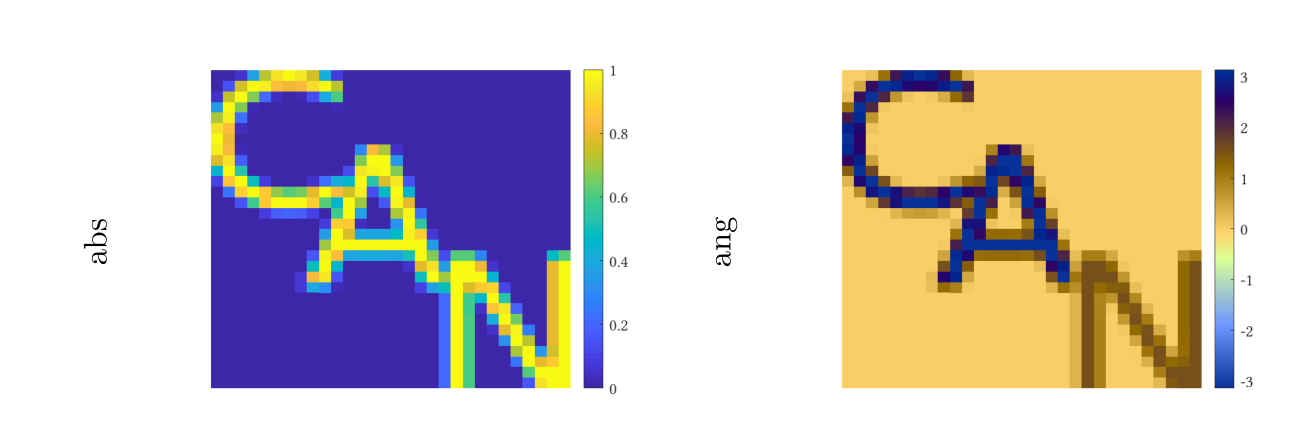}
\caption{Ground truth for amplitude (left) and phase (right, from \cite{Nickel}). \label{phase1}}
\end{figure}

\begin{figure}[t]
\centering
\includegraphics[width=.9\textwidth]{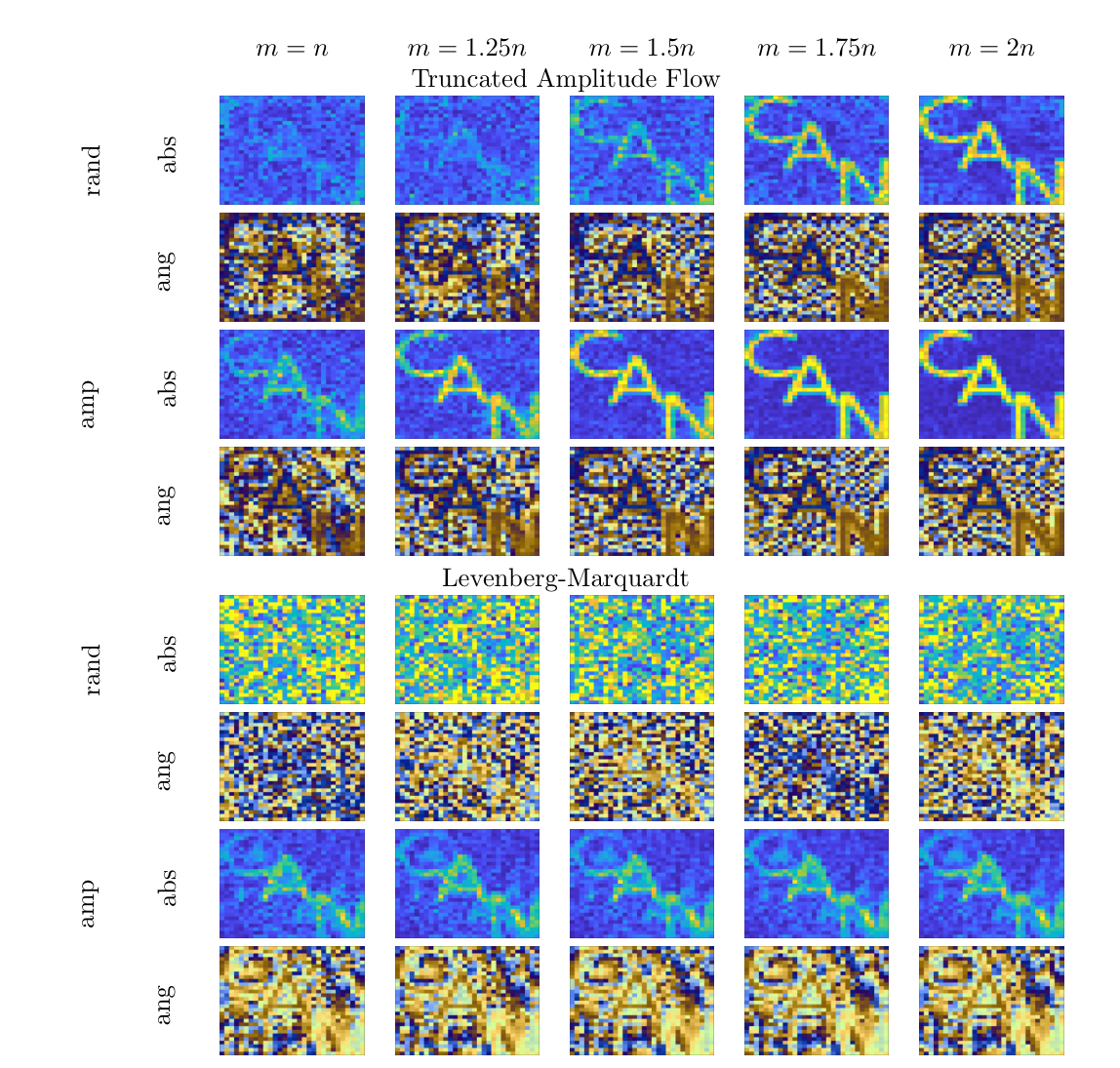}
\caption{Reconstructions of amplitude (left) and phase (right, from \cite{Nickel}) with truncated amplitude flow and Levenberg-Marquardt method with different initializations and sampling rates. \label{phase2}}
\end{figure}

\begin{figure}[t]
	\centering
\includegraphics[width=0.75\textwidth]{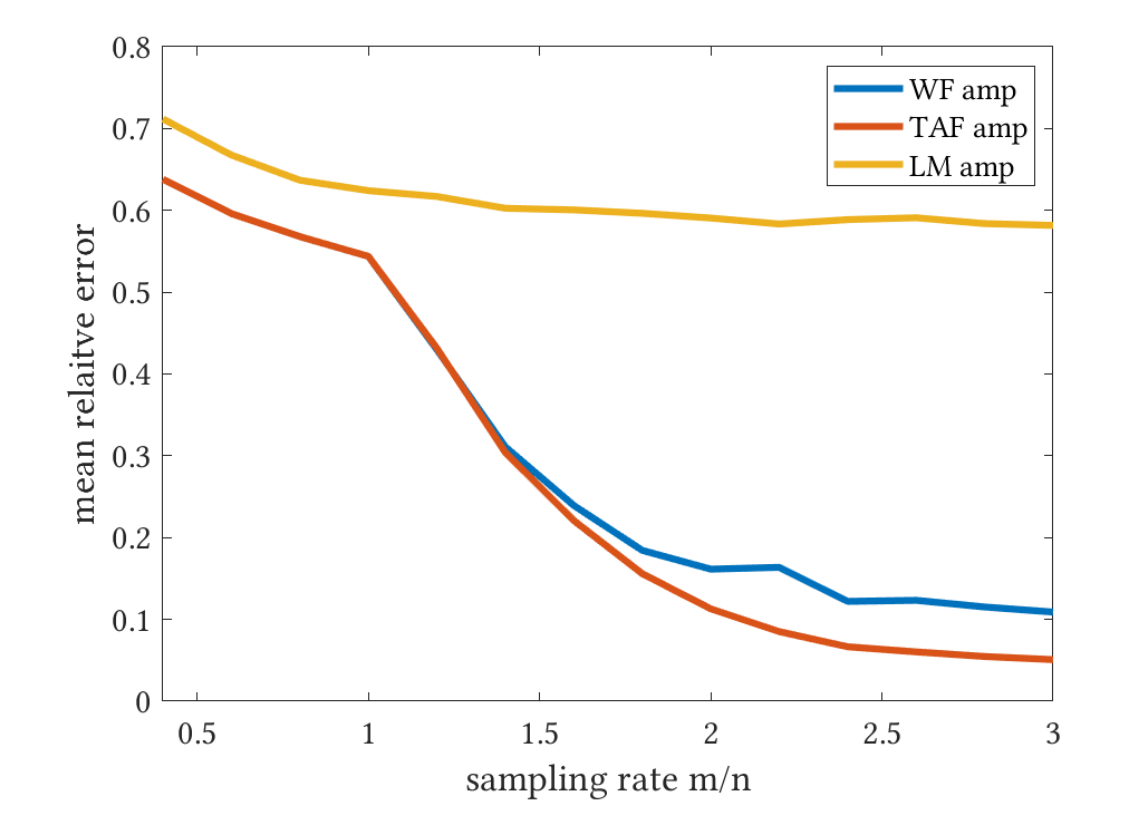}
\caption{Relative mean-square error in the reconstructions of the complex image (from \cite{Nickel}). \label{phase3}}
\end{figure}

Let us finally turn our attention to the reconstruction of the phase in a complex image.
The ground truth for the amplitude and phase of the complex signal are shown in Figure \ref{phase1}.
For brevity we only provide visualizations of reconstructions for the truncated amplitude flow and the Levenberg-Marquardt method in Figure \ref{phase2}, which clearly indicates that the truncated amplitude flow outperforms the Levenberg-Marquardt method, which is the case for all conducted experiments.
Again, a more quantitative evaluation is given in Figure \ref{phase3}, which provides a plot of the relative mean-square error (made invariant to a global phase that is not identifiable) vs. the sampling rate.
We see that the Wirtinger flow can obtain the same reconstruction quality for very low sampling rates, but is outperformed by the truncated amplitude flow at higher sampling rates.

\section{Conclusion}

We have seen that the compressed sensing approach based on single-pixel imaging has great potential to decrease measurement time and effort in THz imaging, but the application in practice depends on several further challenges to be mastered.
First of all appropriate regularization models and numerical algorithms are needed for the image reconstruction problem in order to obtain higher image resolution at reasonable computational times.
Moreover, in several situations it is crucial to consider auto-calibration issues in particular related to the fact that the illuminating beam is difficult to be characterized, and in some cases also diffraction becomes relevant, which effectively yields a phase retrieval problem.
Both effects change the image reconstruction from a problem with a rather simple and incoherent linear forward model to a nonlinear problem with a more coherent forward model, which raises novel computational and also theoretical issues, since the assumptions of the existing compressed sensing theory are not met.

Besides the above mentioned issues several further aspects of practical imaging are foreseen to become relevant for THz single pixel imaging, examples being motion correction for imaging in the field or of moving targets and the reconstruction of multi- or hyper-spectral images, which is a natural motivation in the THz regime.
In the latter case it will become a central question how to combine as few masks as possible in different spectral locations.

\section*{Acknowledgement}
This work has been supported by the German research foundation (DFG) through SPP 1798, Compressed Sensing in Information Processing, projects BU 2327/4-1 and JU 2795/3.
MB acknowledges further support by ERC via Grant EU FP7 ERC Consolidator Grant 615216 LifeInverse.

\bibliographystyle{alpha}
\bibliography{refs}	

\end{document}